\documentclass[leqno,10pt]{article}
\usepackage{theorem}
\usepackage{amssymb}
\usepackage{latexsym}

\theorembodyfont{\itshape}
\newtheorem{Theorem}{Theorem}[section]
\newtheorem{Proposition}{Proposition}[section]

\newtheorem{Lemma}{Lemma}[section]
\newtheorem{Claim}{Claim}[section]
\newtheorem{con}{Conjecture}

\newtheorem{prob}{Problem}
\newtheorem{Corollary}{Corollary}[section]

\newtheorem{s-Theorem}{Theorem}[subsection]
\newtheorem{s-Proposition}{Proposition}[subsection]
\newtheorem{s-Conjecture}{Conjecture}[subsection]
\newtheorem{s-Lemma}{Lemma}[subsection]
\newtheorem{s-Claim}{Claim}[subsection]
\newtheorem{s-Problem}{Problem}[subsection]
\newtheorem{s-Corollary}{Corollary}[subsection]

\theorembodyfont{\rmfamily}
\newtheorem{Definition}{Definition}[section]
\newtheorem{Notation}{Notation}[section]
\newtheorem{Remark}{Remark}[section]
\newtheorem{Example}{Example}[section]
\newtheorem{s-Definition}{Definition}[subsection]
\newtheorem{s-Notation}{Notation}[subsection]
\newtheorem{s-Remark}{Remark}[subsection]
\newtheorem{s-Example}{Example}[subsection]
\newtheorem{definition}{Definition}

\begin{document}
\title{{\bf Effective non-vanishing of global sections of multiple adjoint bundles for polarized $3$-folds}
\thanks{This research was partially supported by the Grant-in-Aid for Young Scientists (B) (No.14740018), The Ministry of Education, Culture, Sports, Science and Technology, Japan, and the Grant-in-Aid for Scientific Research (C)
(No.20540045), Japan Society for the Promotion of Science, Japan.}}
\author{YOSHIAKI FUKUMA}
\date{}
\maketitle
\begin{abstract}
Let $X$ be a smooth complex projective variety of dimension $3$ and let $L$ be an ample line bundle on $X$.
In this paper, we provide a lower bound of $h^{0}(m(K_{X}+L))$ with $\kappa(K_{X}+L)\geq 0$.
In particular, we get the following: (1) if $0\leq \kappa(K_{X}+L)\leq 2$, then $h^{0}(K_{X}+L)>0$ holds. (2) If $\kappa(K_{X}+L)=3$, then $h^{0}(2(K_{X}+L))\geq 3$ holds. Moreover we get a classification of $(X,L)$ with $\kappa(K_{X}+L)=3$ and $h^{0}(2(K_{X}+L))=3$ and $4$.
\end{abstract}

\section{Introduction}
\noindent
\par
Let $X$ be a smooth projective variety of dimension $n$ and let $L$ be an ample (resp. nef and big) line bundle on $X$.
Then the pair $(X,L)$ is called a polarized (resp. quasi-polarized) manifold.
\par
For this $(X,L)$, adjoint bundles $K_{X}+tL$ play important roles 
for investigating this $(X,L)$ 
(for example, see \cite[Chapter 7, 9, and 11]{BeSo-Book}), 
where $K_{X}$ is the canonical line bundle of $X$.
In particular, it is important to know the value of $h^{0}(K_{X}+tL)$.
\par
In \cite[Conjecture 7.2.7]{BeSo-Book}, Beltrametti and Sommese proposed the following conjecture.
\begin{con}\label{Conjecture1}
Let $(X,L)$ be a polarized manifold of dimension $n$.
Assume that $K_{X}+(n-1)L$ is nef.
Then $h^{0}(K_{X}+(n-1)L)>0$.
\end{con}

In \cite[Theorem 2.4]{Fukuma06}, the author proved that Conjecture \ref{Conjecture1}
is true for the case where $\dim X=3$. (See also \cite{BrHo09}.)
Moreover we gave a classification of $(X,L)$ with $h^{0}(K_{X}+2L)=1$ (see \cite[Theorem 2.4]{Fukuma06}).
\par
In general, there is the following conjecture (\cite[Section 4]{Ambro}, \cite[Conjecture 2.1]{Kawamata}).

\begin{con}[Ambro, Kawamata]\label{Conjecture2}
Let $X$ be a complex normal variety, 
$B$ an effective $\mathbb{R}$-divisor on $X$ 
such that the pair $(X,B)$ is KLT, and $D$ a Cartier divisor on $X$.
Assume that $D$ is nef, and that $D-(K_{X}+B)$ is nef and big.
Then $h^{0}(D)>0$.
\end{con}

Here we note that in \cite[Open problems, P.321]{LPS93} Ionescu proposed the same conjecture for the case where $X$ is smooth and $B=0$.

For Conjecture \ref{Conjecture2}, the following results have been obtained.
\begin{itemize}
\item [\rm (\ref{Conjecture2}.a)]
If $\dim X=2$, then Conjecture \ref{Conjecture2} is true (see \cite[Theorem 3.1]{Kawamata}).
\item [\rm (\ref{Conjecture2}.b)]
Let $X$ be a $3$-dimensional projective variety with at most canonical singularities such that $K_{X}$ is nef, and let $D$ be a Cartier divisor such that $D-K_{X}$ is nef and big. Then $h^{0}(D)>0$
 (see \cite[Proposition 4.1]{Kawamata}).
\item [\rm (\ref{Conjecture2}.c)]
Let $(X,L)$ be a polarized manifold of dimension $3$.
Assume that $L^{n}>27$.
Then $h^{0}(K_{X}+L)>0$ if $K_{X}+L$ is nef (see \cite[Th\'eor\`eme 1.8]{Broustet09}).
\item [\rm (\ref{Conjecture2}.d)]
Let $X$ be a $4$-dimensional projective variety with at most Gorenstein canonical singularities. Assume that $D\sim -K_{X}$ is ample. Then $h^{0}(D)>0$
 (see \cite[Theorem 5.2]{Kawamata}).
\item [\rm (\ref{Conjecture2}.e)]
Let $X$ be a smooth projective variety of dimension $3$ with $h^{1}(\mathcal{O}_{X})>0$, and $L$ a nef and big Cartier divisor on $X$ 
such that $K_{X}+L$ is nef.
Then $h^{0}(K_{X}+L)>0$ (see \cite[Theorem 4.2]{1.5}).
\item [\rm (\ref{Conjecture2}.f)]
Let $X$ be a smooth projective variety of dimension $3$ with $\kappa(X)\geq 0$, and $L$ an ample Cartier divisor on $X$.
Then $h^{0}(K_{X}+L)>0$ (see \cite[Theorem 3.2]{Fukuma10}).
\end{itemize}

If $K_{X}+L$ is nef, then by \cite{Shokurov86} there exists a positive integer $m$ such that $h^{0}(m(K_{X}+L))>0$.
More generally if $\kappa(K_{X}+L)\geq 0$, then $h^{0}(m(K_{X}+L))>0$ for some positive integer $m$.
So it is interesting to study the following problem, which was proposed in \cite[Problem 3.2]{Fukuma10}:

\begin{prob}\label{Problem2.8}
For any fixed positive integer $n$, determine the smallest positive integer 
$p$, which depends only on $n$, such that the following {\rm ($*$)} is satisfied:
\begin{itemize}
\item [\rm ($*$)]
$h^{0}(p(K_{X}+L))>0$ for any polarized manifold $(X,L)$ of dimension $n$ 
with $\kappa(K_{X}+L)\geq 0$.
\end{itemize}
\end{prob}
Here we note that
by \cite[Theorem 2.8]{Fukuma10}, 
we see that $p=1$ if $X$ is a curve or surface.
\par
In order to study this problem, in \cite[Problem 5.2]{Fukuma08-3}, we introduced the following:

\begin{definition}
For any fixed positive integer $n$, we set
\begin{eqnarray*}
\mathcal{P}_{n}
&:=&\left \{\ \mbox{\rm $(X,L)$ : polarized manifold}\ |\  \mbox{\rm $\dim X=n$
and
$\kappa(K_{X}+L)\geq 0$}\right\}, \\
\mathcal{M}_{n}
&:=&\left\{\ r\in\mathbb{N}\ |\ h^{0}(r(K_{X}+L))>0\  \mbox{\rm for any $(X,L)\in{\mathcal{P}}_{n}$}\right\},\\
m(n)&:=&
\left\{
\begin{array}{ll}
\mbox{\rm {min}}\ \mathcal{M}_{n}
& \ \ \mbox{\rm if $\mathcal{M}_{n}\neq \emptyset$,} \\
\infty & \ \ \mbox{\rm if $\mathcal{M}_{n}=\emptyset$.}
\end{array}\right. 
\end{eqnarray*}
\end{definition}

In this paper, as the first step, we mainly consider the case where $\dim X=3$.
\par
In \cite[Corollary 5.2]{Fukuma08-3}, we said that $m(3)\leq 2$ holds.
Concretely, in \cite[Theorem 5.4 (2)]{Fukuma08-3}, we proved that if $\kappa(K_{X}+L)=3$, 
then $h^{0}(2(K_{X}+L))\geq 3$.
Moreover in \cite[Theorem 5.4 (1)]{Fukuma08-3}, we announced that in this paper we will prove that $h^{0}(K_{X}+L)>0$ if $0\leq\kappa(K_{X}+L)\leq 2$.
\par
So in this paper, we will prove that $h^{0}(K_{X}+L)>0$ if $n=3$ and $0\leq \kappa(K_{X}+L)\leq 2$.
Moreover, we also study a lower bound of $h^{0}(m(K_{X}+L))$ if $\kappa(K_{X}+L)\geq 0$.
\par
The contents of this paper are the following:
In sections \ref{S2} and \ref{S3}, we will state some definitions and results which will be used later. In particular, in section \ref{S3}, we review the sectional geometric genus.
In section \ref{S4}, we will treat special cases.
If $\kappa(K_{X}+L)=1$ (resp. $2$), then there exists a polarized manifold $(M,A)$ such that $h^{0}(m(K_{X}+L))=h^{0}(m(K_{M}+A))$ for any positive integer $m$ and there exist a fiber space $M\to Y$ such that $Y$ is a normal projective variety of dimension $1$ (resp. $2$), and an ample line bundle $H$ on $Y$ such that $K_{M}+A=f^{*}(H)$. (This $(M,A)$ is called a reduction of $(X,L)$. See Definition \ref{Definition1.7}.)
Hence it is important to consider the following case:
Let $(X,L)$ be a polarized manifold of dimension $n\geq 3$ and let $Y$ be a normal projective variety of dimension $1$ or $2$.
Assume that there exists a fiber space $f:X\to Y$ such that $K_{X}+L=f^{*}(H)$ for some ample line bundle $H$ on $Y$.
In section \ref{S4}, we consider $(X,L)$ like this and we will give a lower bound for $h^{0}(m(K_{X}+L))$.
In particular, we see that $h^{0}(K_{X}+L)>0$ in this case.
\par
In section \ref{S5}, we will study the case where $\dim X=3$.
In particular, we will give a lower bound of $h^{0}(m(K_{X}+L))$ for the following cases:
\begin{itemize}
\item [\rm (a)] $0\leq \kappa(K_{X}+L)\leq 2$ and $m\geq 1$.
\item [\rm (b)] $\kappa(K_{X}+L)=3$ and $m\geq 2$.
\end{itemize}

In particular we get $h^{0}(K_{X}+L)>0$ if $0\leq \kappa(K_{X}+L)\leq 2$ and
$h^{0}(2(K_{X}+L))\geq 3$ if $\kappa(K_{X}+L)=3$ (see also \cite[Theorem 5.4 (2)]{Fukuma08-3}).

Moreover we will also classifiy $(X,L)$ with $\kappa(K_{X}+L)=3$ and $h^{0}(2(K_{X}+L))=3$ or $4$ (see Theorems \ref{EC1} and \ref{EC2}).
\\
\par
In this paper, we shall study mainly a smooth projective variety $X$ 
over the field of complex numbers $\mathbb{C}$.
We will employ the customary notation in algebraic geometry.

\section{Preliminaries}\label{S2}

Here we list up several results which will be used later.

\begin{Definition}\label{Definition1.7}
(i) Let $X$ (resp. $Y$) be an $n$-dimensional projective manifold, and $L$ (resp. $A$) an ample line bundle on $X$ (resp. $Y$).
Then $(X,L)$ is called a  {\it simple blowing up of $(Y,A)$} if there exists a birational morphism $\pi: X\to Y$ such that $\pi$ is a blowing up at a point of $Y$ and $L=\pi^{*}(A)-E$, where $E$ is the $\pi$-exceptional effective reduced divisor.
\\
(ii) Let $X$ (resp. $M$) be an $n$-dimensional projective manifold, and $L$ (resp. $A$) an ample line bundle on $X$ (resp. $M$).
Then we say that $(M,A)$ is a {\it reduction of $(X,L)$} if 
there exists a birational morphism $\mu: X\to M$ such that $\mu$ is a composition of simple blowing ups
and $(M,A)$ is not obtained by a simple blowing up of any polarized manifold.
The map $\mu: X\to M$ is called the {\it reduction map}.
\end{Definition}

\begin{Remark}\label{Remark1.7.1}
Let $(X,L)$ be a polarized manifold and let $(M,A)$ be a reduction of $(X,L)$.
Let $\mu: X\to M$ be the reduction map.
\begin{itemize}
\item [(i)] 
If $(X,L)$ is not obtained by a simple blowing up of another polarized manifold, then $(X,L)$ is a reduction of itself.
\item [(ii)] 
A reduction of $(X,L)$ always exists (see \cite[Chapter II, (11.11)]{Fujita-Book}).
\end{itemize}
\end{Remark}

\begin{Definition}\label{T-D1}
A quasi-polarized surface $(S,L)$ is said to be {\it $L$-minimal} if $LE>0$ for every $(-1)$-curve $E$ on $S$.
\end{Definition}

\begin{Lemma}\label{Lemma B}
Let $X$ be a complete normal variety of dimension $n$, 
and let $D_{1}$ and $D_{2}$ be effective Cartier divisors on $X$.
Then $h^{0}(D_{1}+D_{2})\geq h^{0}(D_{1})+h^{0}(D_{2})-1$.
\end{Lemma}
\noindent{\em Proof.}
See \cite[Lemma 1.10]{Fukuma3} or \cite[15.6.2 Lemma]{Kollar95}. $\Box$

\begin{Proposition}\label{GHIT}
Let $X$ be a projective variety of dimension $n$ and let $D_{i}$ be $\mathbb{Q}$-Cartier divisors on $X$ for $0\leq i\leq k$.
Assume that $n\geq 2$ and that $D_{i}$ is nef for every integer $i$ with $1\leq i\leq k$.
If $n_{1}+\cdots +n_{k}=n-1$ and $n_{1}\geq 1$, then we have
$$(D_{0}D_{1}^{n_{1}}\cdots D_{k}^{n_{k}})^{2}\geq 
(D_{0}^{2}D_{1}^{n_{1}-1}\cdots D_{k}^{n_{k}})(D_{1}^{n_{1}+1}\cdots D_{k}^{n_{k}}).$$
\end{Proposition}
\noindent{\em Proof.}
See \cite[Proposition 2.5.1]{BeSo-Book}. $\Box$

\begin{Proposition}\label{2P-T1}
Let $X$ be a normal projective surface and let $\pi: S\to X$ be a resolution of singularities of $X$.
Then $\chi(\mathcal{O}_{S})+h^{0}(R^{1}\pi_{*}(\mathcal{O}_{S}))=\chi(\mathcal{O}_{X})$.
In particular $\chi(\mathcal{O}_{S})\leq \chi(\mathcal{O}_{X})$ holds.
\end{Proposition}
\noindent
{\em Proof.}
By using Leray's spectral sequence for $\pi^{*}(\mathcal{O}_{X})$,
we have 
$$\chi(\pi^{*}\mathcal{O}_{X})=\sum_{q\geq 0}(-1)^{q}\chi(R^{q}\pi_{*}(\pi^{*}\mathcal{O}_{X})).$$
Since $R^{q}\pi_{*}(\pi^{*}\mathcal{O}_{X})\cong R^{q}\pi_{*}(\mathcal{O}_{S})$
and $R^{q}\pi_{*}(\mathcal{O}_{S})=0$ for every integer $q$ with $q\geq 2$,
we have 
$$\chi(\pi^{*}\mathcal{O}_{X})=\chi(\pi_{*}(\mathcal{O}_{S}))-\chi(R^{1}\pi_{*}(\mathcal{O}_{S})).$$
Here we also note that $\pi_{*}(\mathcal{O}_{S})=\mathcal{O}_{X}$ because $\pi$ is birational and $X$ is normal (see \cite[Corollary 11.4 in Chapter III]{Hartshorne}).
Moreover $\chi(R^{1}\pi_{*}(\mathcal{O}_{S}))=h^{0}(R^{1}\pi_{*}(\mathcal{O}_{S}))$ because $\dim\mbox{Supp}(R^{1}\pi_{*}(\mathcal{O}_{S}))\leq 0$.
Therefore since $\mathcal{O}_{S}=\pi^{*}(\mathcal{O}_{X})$, we get the assertion. $\Box$

\begin{Lemma}\label{4L-T4}
Let $X$ be a smooth projective variety of dimension $n$ and let $Y$ be a normal projective variety of dimension $m$ with $n>m\geq 1$.
Assume that $q(X)=q(Y)$ and there exists a fiber space $f: X\to Y$,
that is, $f$ is a surjective morphism with connected fibers.
Then for any resolution of singularities of $Y$, $\pi: Z\to Y$, we have $q(Z)=q(Y)$.
In particular, if $q(Y)\geq 1$, then the Albanese map of $Y$ can be defined.
\end{Lemma}
\noindent
{\em Proof.} By assumption, there exist smooth projective varieties $X_{1}$ and $Y_{1}$, birational morphisms $\mu_{1}: X_{1}\to X$ and $\nu_{1}: Y_{1}\to Y$,
and a fiber space $f_{1}:X_{1}\to Y_{1}$ 
such that $f\circ \mu_{1}=\nu_{1}\circ f_{1}$.
Here we note that $q(X)=q(X_{1})$ and $q(X_{1})\geq q(Y_{1})$.
Moreover we have $q(Y_{1})\geq q(Y)$ holds.
Hence we get $q(Y_{1})\geq q(Y)=q(X)=q(X_{1})\geq q(Y_{1})$ and we have $q(Y_{1})=q(Y)$.
On the other hand let $Z$ be any resolution of singularities of $Y$.
Then $q(Z)=q(Y_{1})$ because $Z$ is birationally equivalent to $Y_{1}$.
In particular, by \cite[(0.3.3) Lemma]{12} or \cite[Lemma 2.4.1 and Remark 2.4.2]{BeSo-Book}, the Albanese map of $Y$ can be defined.
Hence we get the assertion of Lemma \ref{4L-T4}. $\Box$

\section{Review on the sectional geometric genus}\label{S3}
In this section, we review the definition and some properties of the sectional
geometric genus of polarized manifolds, which will be used later.

\begin{Notation}\label{Notation1.1}
Let $X$ be a projective variety of dimension $n$ 
and let $L$ be a line bundle on $X$.
Let $\chi(tL)$ be the Euler-Poincar\'e characteristic of $tL$, where $t$ is an indeterminate.
Then we put
$$\chi(tL)=\sum_{j=0}^{n}\chi_{j}(X,L){t+j-1\choose j}.$$
\end{Notation}

\begin{Definition}\label{Definition1.2}
Let $X$ be a projective variety of dimension $n$ 
and let $L$ be a line bundle on $X$.
Then for every integer $i$ with $0\leq i\leq n$, 
the {\it $i$-th sectional $H$-arithmetic genus $\chi_{i}^{H}(X,L)$} and 
the {\it $i$-th sectional geometric genus $g_{i}(X,L)$ of $(X,L)$} are defined by the following:
\begin{eqnarray*}
\chi_{i}^{H}(X,L)&:=&\chi_{n-i}(X,L),\\
g_{i}(X,L)
&:=&(-1)^{i}(\chi_{i}^{H}(X,L)-\chi(\mathcal{O}_{X}))
+\sum_{j=0}^{n-i}(-1)^{n-i-j}h^{n-j}(\mathcal{O}_{X}).
\end{eqnarray*}
\end{Definition}

\begin{Remark}\label{Remark1.2.1}
\begin{itemize}
\item [(1)] Since $\chi_{n-i}(X,L)\in \mathbb{Z}$, we see that $\chi_{i}^{H}(X,L)$ and $g_{i}(X,L)$ are integers by definition.
\item [(2)] 
If $i=0$, then $\chi_{0}^{H}(X,L)$ and $g_{0}(X,L)$ are equal to the degree of $(X,L)$.
\item [(3)] 
If $i=1$, then $g_{1}(X,L)$ is equal to the sectional genus $g(X,L)$ of $(X,L)$.
\item [(4)] 
If $i=n$, then $\chi_{n}^{H}(X,L)=\chi(\mathcal{O}_{X})$ and $g_{n}(X,L)=h^{n}(\mathcal{O}_{X})$.
\end{itemize}
\end{Remark}

\begin{Theorem}\label{Theorem1.3}
Let $(X,L)$ be a quasi-polarized manifold with $\dim X=n$.
For every integer $i$ with $0\leq i\leq n-1$, we have
$$g_{i}(X,L)=\sum_{j=0}^{n-i-1}(-1)^{j}{n-i\choose j}h^{0}(K_{X}+(n-i-j)L)
+\sum_{k=0}^{n-i}(-1)^{n-i-k}h^{n-k}(\mathcal{O}_{X}).$$
\end{Theorem}
\noindent{\em Proof.} 
See \cite[Theorem 2.3]{Fukuma3}. $\Box$
\\
\par The following theorem will be often used later.
\begin{Theorem}\label{Theorem1.4}
Let $(X,L)$ be a polarized $3$-fold.
Assume that $\kappa(K_{X}+L)\geq 0$.
Then $g_{2}(X,L)\geq h^{1}(\mathcal{O}_{X})$.
\end{Theorem}
\noindent
{\em Proof.}
See \cite[Theorem 3.3.1 (2)]{Fukuma05-2}. $\Box$

\begin{Notation}\label{B1}
Let $X$ be a projective variety of dimension $n$, let $i$ be an integer with $0\leq i\leq n-1$, and let $L_{1},\dots , L_{n-i}$ be line bundles on $X$.
Then $\chi(L_{1}^{t_{1}}\otimes\cdots\otimes L_{n-i}^{t_{n-i}})$ is a polynomial in $t_{1}, \dots ,t_{n-i}$ of total degree at most $n$.
So we can write $\chi(L_{1}^{t_{1}}\otimes\cdots
\otimes L_{n-i}^{t_{n-i}})$ uniquely as follows.
\begin{eqnarray*}
&&\chi(L_{1}^{t_{1}}\otimes\cdots\otimes L_{n-i}^{t_{n-i}}) \\
&&=\sum_{p=0}^{n}\sum
_{\stackrel{p_{1}\geq 0,\dots , p_{n-i}\geq 0}
{p_{1}+\cdots +p_{n-i}=p}}
\chi_{p_{1},\dots , p_{n-i}}(L_{1},\dots ,L_{n-i})
{t_{1}+p_{1}-1\choose p_{1}}\dots{t_{n-i}+p_{n-i}-1\choose p_{n-i}}.
\end{eqnarray*}
\end{Notation}

\begin{Definition}\label{B2}(\cite[Definition 2.1 and Remark 2.2 (2)]{Fukuma11})
Let $X$ be a projective variety of dimension $n$, let $i$ be an integer with $0\leq i\leq n$, and let $L_{1},\dots , L_{n-i}$ be line bundles on $X$.
\\
(1) The {\it $i$-th sectional $H$-arithmetic genus $\chi_{i}^{H}(X,L_{1},\dots , L_{n-i})$} is defined by the following:

\[
\chi_{i}^{H}(X,L_{1},\dots , L_{n-i})=
\left\{
\begin{array}{ll}
\chi_{\underbrace{1, \dots , 1}_{n-i}}(L_{1},\dots , L_{n-i}) 
& \mbox{if $0\leq i\leq n-1$,} \\
\chi(\mathcal{O}_{X}) & \mbox{if $i=n$.}
\end{array}\right. \] 
\\
(2) The {\it $i$-th sectional geometric genus $g_{i}(X,L_{1},\dots , L_{n-i})$} is defined by the following:
\begin{eqnarray*}
g_{i}(X,L_{1},\dots , L_{n-i})
&=&(-1)^{i}(\chi_{i}^{H}(X,L_{1},\dots , L_{n-i})-\chi({\mathcal{O}}_{X})) \\
&&\ \ \ +\sum_{j=0}^{n-i}(-1)^{n-i-j}h^{n-j}({\mathcal{O}}_{X}).
\end{eqnarray*}
\end{Definition}

\begin{Remark}\label{B6}
(1) Let $X$ be a projective variety of dimension $n$ and 
let $L$ be a line bundle on $X$.
Let $i$ be an integer with $0\leq i\leq n-1$.  
Then 
$$\chi_{i}^{H}(X,L,\dots , L)=\chi_{i}^{H}(X,L)$$
and 
$$g_{i}(X,L,\dots , L)=g_{i}(X,L).$$ 
(See \cite[Corollary 2.1]{Fukuma11}.)
\\
(2)
Let $X$ be a smooth projective variety of dimension $n$, and let $L_{1},\dots , L_{n-1}$ be line bundles on $X$.
Then 
$$g_{1}(X,L_{1}, \dots, L_{n-1})=1+\frac{1}{2}\left(K_{X}+\sum_{j=1}^{n-1}L_{j}\right)L_{1}\cdots L_{n-1}.$$
(See \cite[Corollary 2.7]{Fukuma11} or \cite[Proposition 6.1.1]{Fukuma08-4}.)
\end{Remark}

\begin{Theorem}\label{B18}
Let $i$ be an integer with $1\leq i\leq n$.
Let $A, B, L_{1}, \cdots , L_{n-i-1}$ be line bundles on $X$.
Then
\begin{eqnarray*}
&&\chi_{i}^{H}(X,A+B,L_{1},\cdots , L_{n-i-1})\\
&&=\chi_{i}^{H}(X,A,L_{1},\cdots , L_{n-i-1})+\chi_{i}^{H}(X,B,L_{1},\cdots , L_{n-i-1}) \\
&&\ \ \ -\chi_{i-1}^{H}(X,A,B,L_{1},\cdots , L_{n-i-1})
\end{eqnarray*}
\begin{eqnarray*}
&&g_{i}(X,A+B,L_{1},\cdots , L_{n-i-1}) \\
&&=g_{i}(X,A,L_{1},\cdots , L_{n-i-1})+g_{i}(X,B,L_{1},\cdots , L_{n-i-1}) \\
&&\ \ \ +g_{i-1}(X,A,B,L_{1},\cdots , L_{n-i-1})-h^{i-1}(\mathcal{O}_{X}).
\end{eqnarray*}
\end{Theorem}
\noindent{\em Proof.} See \cite[Corollary 2.4]{Fukuma11}. $\Box$

\begin{Proposition}\label{SKB}
Let $X$ be a smooth projective variety with $\dim X=n\geq 2$,
let $L_{1},\cdots ,L_{m}$ be nef and big line bundles on $X$ and let $L$ be a nef line bundle, where $m\geq 1$.
Then
\begin{eqnarray*}
&&h^{0}(K_{X}+L_{1}+\cdots +L_{m}+L)-h^{0}(K_{X}+L_{1}+\cdots +L_{m})\\
&&=\sum_{s=0}^{n-1}\sum_{(k_{1},\cdots,k_{n-s-1})\in A_{n-s-1}^{m}}
g_{s}(X, L_{k_{1}},\cdots, L_{k_{n-s-1}},L) \\
&&\ \ \ -\sum_{s=0}^{n-2}{m-1\choose n-s-2}h^{s}(\mathcal{O}_{X}).
\end{eqnarray*}
Here $A_{t}^{p}:=\left\{ (k_{1},\cdots , k_{t})\ |\ k_{l}\in \{ 1, \cdots , p\}, k_{i}<k_{j} \ \mbox{\rm if $i<j$}\right\}$, and we set
\[
\sum_{(k_{1},\cdots,k_{n-s-1})\in A_{n-s-1}^{m}}
g_{s}(X, L_{k_{1}},\cdots, L_{k_{n-s-1}},L)
=\left\{
\begin{array}{lc}
0 & \mbox{if $n-s-1>m$,} \\
g_{n-1}(X,L) & \mbox{if $s=n-1$.}
\end{array} \right. \]
\end{Proposition}
\noindent{\em Proof.}
See \cite[Theorem 5.1]{Fukuma08-3}. $\Box$

\section{Special cases}\label{S4}

In this section, we will investigate the dimension of adjoint linear system for special cases.
First we prove the following.

\begin{Theorem}\label{3-1T1}
Let $(X,L)$ be a polarized manifold of dimension $n\geq 2$ and let $C$ be a smooth projective curve.
Assume that there exists a fiber space $f:X\to C$ such that $K_{X}+L=f^{*}(H)$ for some ample line bundle $H$ on $C$.
Then for every positive integer $m$ 
\[
h^{0}(m(K_{X}+L))\geq
\left\{
\begin{array}{lc}
(m-1)(g(C)-1)+mg(C) & \mbox{if $g(C)\geq 1$,} \\
m+1 & \mbox{if $g(C)=0$.}
\end{array} \right. 
\]
In particular $h^{0}(K_{X}+L)>0$ holds.
\end{Theorem}
\noindent{\em Proof.}
In this case
\begin{eqnarray*}
h^{0}(m(K_{X}+L))
&=&h^{0}(f^{*}(mH)) \\
&=&h^{0}(mH) \\
&=&h^{1}(mH)+\deg(mH)+(1-g(C)).
\end{eqnarray*}
On the other hand, by \cite[Lemma 1.13]{Fukuma3}, we have $\deg H\geq 2g(C)-1$.
Hence if $g(C)\geq 1$, then 
\begin{eqnarray*}
h^{0}(mH)&\geq& m(2g(C)-1)+1-g(C)\\
&=&(2m-1)g(C)-(m-1) \\
&=&(m-1)(g(C)-1)+mg(C).
\end{eqnarray*}
If $g(C)=0$, then $h^{1}(mH)=0$ and 
$h^{0}(mH)=\deg(mH)+1\geq m+1$.
Therefore
\[
h^{0}(m(K_{X}+L))\geq
\left\{
\begin{array}{lc}
(m-1)(g(C)-1)+mg(C) & \mbox{if $g(C)\geq 1$,} \\
m+1 & \mbox{if $g(C)=0$.}
\end{array} \right. 
\]
This completes the proof. $\Box$

\begin{Corollary}\label{3-1T1C}
Let $(X,L)$ be a polarized manifold of dimension $n\geq 2$ and let $C$ be a smooth projective curve.
Assume that there exists a fiber space $f:X\to C$ such that $K_{X}+L=f^{*}(H)$ for some ample line bundle $H$ on $C$.
Then for every positive integer $m$ 
\[
h^{0}(m(K_{X}+L))\geq
\left\{
\begin{array}{lc}
m & \mbox{if $g(C)\geq 1$,} \\
m+1 & \mbox{if $g(C)=0$.}
\end{array} \right. 
\]
\end{Corollary}

\begin{Theorem}\label{3-1T3}
Let $(X,L)$ be a polarized manifold of dimension $n\geq 2$ and let $C$ be a smooth projective curve.
Assume that there exists a fiber space $f:X\to C$ such that $K_{X}+L=f^{*}(H)$ for some ample line bundle $H$ on $C$.
\begin{itemize}
\item [\rm (1)]
If $g(C)\geq 1$ and $h^{0}(m(K_{X}+L))=m$ for some positive integer $m$, then 
$g(C)=1$ and $\deg H=1$.
\item [\rm (2)] 
If $g(C)=0$ and $h^{0}(m(K_{X}+L))=m+1$ for some positive integer $m$, then 
$(C,H)\cong (\mathbb{P}^{1}, \mathcal{O}_{\mathbb{P}^{1}}(1))$.
\end{itemize}
\end{Theorem}
\noindent{\em Proof.}
(2.1) Assume that $g(C)\geq 1$ and $h^{0}(m(K_{X}+L))=m$.
Then by the proof of Theorem \ref{3-1T1} we have $g(C)=1$ and $\deg H=1$.
\\
(2.2) Assume that $g(C)=0$ and $h^{0}(m(K_{X}+L))=m+1$.
Then the proof of Theorem \ref{3-1T1} implies that $\deg H=1$, that is, $H=\mathcal{O}_{\mathbb{P}^{1}}(1)$.
Therefore $(C,H)\cong (\mathbb{P}^{1},\mathcal{O}_{\mathbb{P}^{1}}(1))$.
So we get the assertion. $\Box$
\\
\par
Next we consider the following case.

\begin{Theorem}\label{3-1T2}
Let $(X,L)$ be a polarized manifold of dimension $n\geq 3$ and let $Y$ be a normal projective surface.
Assume that there exists a fiber space $f:X\to Y$ such that $K_{X}+L=f^{*}(H)$ for some ample line bundle $H$ on $Y$.
Then for every positive integer $m$
\[
h^{0}(m(K_{X}+L))\geq
\left\{
\begin{array}{lc}
{m+1\choose 2}-(m-1)\chi(\mathcal{O}_{Y}) & \mbox{if $\chi(\mathcal{O}_{Y})\leq 0$,} \\
{m\choose 2}+\chi(\mathcal{O}_{Y}) & \mbox{if $\chi(\mathcal{O}_{Y})>0$.}
\end{array} \right. 
\]
In particular $h^{0}(K_{X}+L)>0$ holds.
\end{Theorem}
\noindent{\em Proof.}
In this case $h^{0}(m(K_{X}+L))=h^{0}(mH)$.
Here we note the following.
\begin{Claim}\label{CL1}
$h^{i}(mH)=0$ for $i=1, 2$.
\end{Claim}
\noindent
{\em Proof.}
Since $f^{*}(mH)-K_{X}=(m-1)K_{X}+mL=(m-1)(K_{X}+L)+L$ is ample, 
we have $R^{i}f_{*}(f^{*}(mH))=0$ for every $i>0$ 
by \cite[Theorem 1.7]{Fukuma3}.
Hence by \cite[Exsercise 8.1 page 252 in Chapter III]{Hartshorne} we have
$h^{i}(f^{*}(mH))=h^{i}(f_{*}f^{*}(mH))=h^{i}(mH)$.
Therefore for every $i>0$
\begin{eqnarray*}
h^{i}(mH)&=&h^{i}(f^{*}(mH)) \\
&=&h^{i}(m(K_{X}+L)) \\
&=&h^{i}(K_{X}+(m-1)(K_{X}+L)+L) \\
&=&0.
\end{eqnarray*}
This completes the proof of Claim \ref{CL1}. $\Box$
\\
\par
By Claim \ref{CL1}, we have $h^{0}(m(K_{X}+L))=h^{0}(mH)=\chi(mH)$.
Here we use Notation \ref{Notation1.1}.
Then $\chi_{0}(Y,H)=\chi(\mathcal{O}_{Y})$, $\chi_{1}(Y,H)=1-g(Y,H)$ and
$\chi_{2}(Y,H)=H^{2}$, where $g(Y,H)$ denotes the sectional genus of $(Y,H)$.
Let $\delta: S\to Y$ be a minimal resolution of $Y$.
Then there exist a smooth projective variety $X_{1}$, 
a birational morphism $\mu_{1}:X_{1}\to X$ and 
a fiber space $f_{1}:X_{1}\to S$ such that $f\circ \mu_{1}=\delta\circ f_{1}$.
\\
\\
(I) The case where $\chi(\mathcal{O}_{Y})\leq 0$.
\\
Then
\begin{eqnarray}
\chi(mH)-m\chi(H)
&=&\sum_{j=0}^{2}\chi_{j}(Y,H){m+j-1\choose j}-m\sum_{j=0}^{2}\chi_{j}(Y,H)
\label{T-EQ1}\\
&=&-(m-1)\chi(\mathcal{O}_{Y})+\left({m+1\choose 2}-m\right)H^{2}\nonumber\\
&\geq&{m+1\choose 2}-m-(m-1)\chi(\mathcal{O}_{Y}) \nonumber\\
&=&{m\choose 2}-(m-1)\chi(\mathcal{O}_{Y}). \nonumber
\end{eqnarray}
Therefore $\chi(mH)\geq m\chi(H)+{m\choose 2}-(m-1)\chi(\mathcal{O}_{Y})=mh^{0}(H)+{m\choose 2}-(m-1)\chi(\mathcal{O}_{Y})$.
\\
\par
Next we prove the following claim.

\begin{Claim}\label{CL3}
$h^{0}(H)>0$.
\end{Claim}
\noindent{\em Proof.}
Since $\chi(\mathcal{O}_{Y})\leq 0$ in this case, we see that $h^{1}(\mathcal{O}_{Y})>0$. Because $h^{1}(\mathcal{O}_{X})=h^{1}(\mathcal{O}_{Y})$ in this case, by Lemma \ref{4L-T4} we see that $Y$ has the Albanese map.
Let $\alpha: Y\to \mbox{Alb}(Y)$ be the Albanese map of $Y$
and let $h:=\alpha\circ f$.
Here we note that $\dim h(X)=1$ or $2$.
\\
(a) First we consider the case where $\dim h(X)=2$.
By \cite[Corollary 10.7 in Chapter III]{Hartshorne}
any general fiber $F_{h}$ of $h$ can be written as follows:
$F_{h}=\cup_{i=1}^{r}F_{i}$, where $F_{i}$ is a smooth projective variety 
of dimension $n-2$.
We note that $F_{i}$ is a fiber of $f$ for every $i$.
Since $(K_{X}+L)|_{F_{i}}=f^{*}(H)|_{F_{i}}\cong \mathcal{O}_{F_{i}}$, 
we have
$$
h^{0}((K_{X}+L)|_{F_{h}})
=\sum_{i=1}^{r}h^{0}(K_{F_{i}}+L_{F_{i}})
=\sum_{i=1}^{r}h^{0}(\mathcal{O}_{F_{i}})
>0.
$$
By \cite[Lemma 4.1]{1.5} we have $h^{0}(H)=h^{0}(K_{X}+L)>0$.
\\
(b) Next we consider the case where $\dim h(X)=1$.
Then we note that $h$ has connected fibers.
Let $F_{h}$ (resp. $F_{\alpha}$) be a general fiber of $h$ (resp. $\alpha$).
Then $f|_{F_{h}}:F_{h}\to F_{\alpha}$ is a fiber space such that
$K_{F_{h}}+L_{F_{h}}=f^{*}(H)|_{F_{h}}=(f|_{F_{h}})^{*}(H|_{F_{\alpha}})$.
Here we note that $F_{h}$ and $F_{\alpha}$ are smooth projective varieties.
Since $H$ is ample, so is $H_{F_{\alpha}}$ on $F_{\alpha}$.
Since $\dim F_{\alpha}=1$, by Theorem \ref{3-1T1}
we have $h^{0}(K_{F_{h}}+L_{F_{h}})>0$.
Therefore by \cite[Lemma 4.1]{1.5} we get $h^{0}(H)=h^{0}(K_{X}+L)>0$.
This completes the proof. $\Box$
\\
\par
Claim \ref{CL3} implies that by (\ref{T-EQ1})
\begin{eqnarray*}
\chi(mH)&\geq& mh^{0}(H)+{m\choose 2}-(m-1)\chi(\mathcal{O}_{Y}) \\
&\geq& m+{m\choose 2}-(m-1)\chi(\mathcal{O}_{Y}) \\ 
&\geq& {m+1\choose 2}-(m-1)\chi(\mathcal{O}_{Y}).
\end{eqnarray*}
\noindent
\\
(II) Next we consider the case where $\chi(\mathcal{O}_{Y})>0$.
First we prove the following lemma.
\begin{Lemma}\label{L1}
$\chi_{1}(Y,H)+\chi_{2}(Y,H)\geq 0$.
\end{Lemma}
\noindent{\em Proof.}
First we note that
$K_{X_{1}}+\mu_{1}^{*}(L)\geq \mu_{1}^{*}(K_{X}+L)=\mu_{1}^{*}f^{*}(H)=f_{1}^{*}\delta^{*}(H)$.
Hence
for a general fiber $F_{1}$ of $f_{1}$, we have $0<h^{0}((K_{X_{1}}+\mu_{1}^{*}(L))|_{F_{1}})=h^{0}(K_{F_{1}}+(\mu_{1}^{*}(L))_{F_{1}})$.
Hence we have $(f_{1})_{*}(K_{X_{1}/S}+\mu_{1}^{*}(L))\neq 0$.
By Hironaka's theory there exist a smooth projective variety $X_{2}$ and a birational morphism $\mu_{2}:X_{2}\to X_{1}$ such that 
$$\mu_{2}^{*}f_{1}^{*}((f_{1})_{*}(K_{X_{1}/S}+\mu_{1}^{*}(L)))\to \mu_{2}^{*}(K_{X_{1}/S}+\mu_{1}^{*}(L)-D)-E_{2}$$
is surjective, where $D$ is an effective divisor on $X_{1}$ and $E_{2}$ is a $\mu_{2}$-exceptional effective divisor on $X_{2}$.
Since $(f_{1})_{*}(K_{X_{1}/S}+\mu_{1}^{*}(L))$ is weakly positive (\cite[Theorem A$^{\prime}$ in Appendix]{Fukuma1.5}),
we see that $\mu_{2}^{*}(K_{X_{1}/S}+\mu_{1}^{*}(L)-D)-E_{2}$ is pseudo effective (see the proof of (1) in \cite[Remark 1.3.2]{Fukuma1.5}).
Here we note that for every positive integer $p$ we have
$$0\leq (\mu_{2}^{*}(K_{X_{1}/S}+\mu_{1}^{*}(L)-D)-E_{2})\mu_{2}^{*}f_{1}^{*}\delta^{*}(H)(\mu_{2}^{*}\mu_{1}^{*}(pL))^{n-2}$$ 
because $H$ is ample.
On the other hand
\begin{eqnarray*}
&&(\mu_{2}^{*}(K_{X_{1}/S}+\mu_{1}^{*}(L)-D)-E_{2})\mu_{2}^{*}f_{1}^{*}\delta^{*}(H)(\mu_{2}^{*}\mu_{1}^{*}(pL))^{n-2}\\
&&=(K_{X_{1}/S}+\mu_{1}^{*}(L)-D)(f_{1}^{*}\delta^{*}(H))(\mu_{1}^{*}(pL))^{n-2}\\
&&\leq (K_{X_{1}/S}+\mu_{1}^{*}(L))(f_{1}^{*}\delta^{*}(H))(\mu_{1}^{*}(pL))^{n-2}.
\end{eqnarray*}
Since $K_{X_{1}}=\mu_{1}^{*}K_{X}+E_{1}$, where $E_{1}$ is a $\mu_{1}$-exceptional effective divisor on $X_{1}$,
we have
\begin{eqnarray*}
&&(K_{X_{1}/S}+\mu_{1}^{*}(L))(f_{1}^{*}\delta^{*}(H))(\mu_{1}^{*}(pL))^{n-2}\\
&&=(\mu_{1}^{*}(K_{X}+L)-f_{1}^{*}(K_{S})+E_{1})(f_{1}^{*}\delta^{*}(H))(\mu_{1}^{*}(pL))^{n-2} \\
&&=(f_{1}^{*}(\delta^{*}(H)-K_{S})+E_{1})(\mu_{1}^{*}f^{*}(H))(\mu_{1}^{*}(pL))^{n-2} \\
&&=f_{1}^{*}(\delta^{*}(H)-K_{S})(\mu_{1}^{*}f^{*}(H))(\mu_{1}^{*}(pL))^{n-2}\\
&&=f_{1}^{*}(\delta^{*}(H)-K_{S})(f_{1}^{*}\delta^{*}(H))(\mu_{1}^{*}(pL))^{n-2}.
\end{eqnarray*}
Here we take $p$ as $\mbox{Bs}|\mu_{1}^{*}(pL)|=\emptyset$.
Then there exist $(n-2)$-general members $H_{1}, \dots , H_{n-2}$ in $|\mu_{1}^{*}(pL)|$ such that $H_{1}\cap \dots \cap H_{n-2}$ is a smooth projective surface $S_{1}$.
Then $f_{1}|_{S}:S_{1}\to S$ is a surjective morphism and we have
\begin{eqnarray*}
&&f_{1}^{*}(\delta^{*}(H)-K_{S})(f_{1}^{*}\delta^{*}(H))(\mu_{1}^{*}(pL))^{n-2}\\
&&=f_{1}^{*}(\delta^{*}(H)-K_{S})f_{1}^{*}(\delta^{*}(H))S_{1}\\
&&=(\deg f_{1}|_{S_{1}})(\delta^{*}(H)-K_{S})\delta^{*}(H).
\end{eqnarray*}
On the other hand, since $\chi_{2}(Y,H)=\chi_{2}(S,\delta^{*}(H))$ and $\chi_{1}(Y,H)=\chi_{1}(S,\delta^{*}(H))$,
we have
$(\delta^{*}(H)-K_{S})\delta^{*}(H)=2(\chi_{1}(S,\delta^{*}(H))+\chi_{2}(S,\delta^{*}(H)))=2(\chi_{1}(Y,H)+\chi_{2}(Y,H))$.
Hence we get the assertion. $\Box$
\\
\par
Therefore we get
\begin{eqnarray*}
h^{0}(mH)=\chi(mH)
&=& \chi_{0}(Y,H)+\chi_{1}(Y,H)m+\chi_{2}(Y,H){m+1\choose 2}\\
&=& \chi(\mathcal{O}_{Y})+m(\chi_{1}(Y,H)+\chi_{2}(Y,H))+\left({m+1\choose 2}-m\right)\chi_{2}(Y,H)\\
&\geq& \chi(\mathcal{O}_{Y})+{m\choose 2}.
\end{eqnarray*}
Therefore 
$$h^{0}(m(K_{X}+L))\geq {m\choose 2}+\chi(\mathcal{O}_{Y}).$$
This completes the proof. $\Box$

\begin{Corollary}\label{3-1T2C}
Let $(X,L)$ be a polarized manifold of dimension $n\geq 3$ and let $Y$ be a normal projective surface.
Assume that there exists a fiber space $f:X\to Y$ such that $K_{X}+L=f^{*}(H)$ for some ample line bundle $H$ on $Y$.
Then for every positive integer $m$ 
\[
h^{0}(m(K_{X}+L))\geq
\left\{
\begin{array}{lc}
{m+1\choose 2} & \mbox{if $\chi(\mathcal{O}_{Y})\leq 0$,} \\
{m\choose 2}+1 & \mbox{if $\chi(\mathcal{O}_{Y})>0$.}
\end{array} \right. 
\]
\end{Corollary}

\begin{Theorem}\label{3-1T4}
Let $(X,L)$ be a polarized manifold of dimension $n\geq 3$ and let $Y$ be a normal projective surface.
Assume that there exists a fiber space $f:X\to Y$ such that $K_{X}+L=f^{*}(H)$ for some ample line bundle $H$ on $Y$.
\begin{itemize}
\item [\rm (1)]
If $\chi(\mathcal{O}_{Y})\leq 0$ and $h^{0}(m(K_{X}+L))={m+1\choose 2}$ for some positive integer $m\geq 2$, then 
$Y$ is smooth and $(Y,H)$ is a scroll over a smooth elliptic curve $C$ such that $H^{2}=1$.
\item [\rm (2)] 
If $\chi(\mathcal{O}_{Y})>0$ and $h^{0}(m(K_{X}+L))={m\choose 2}+1$ for some positive integer $m\geq 2$, then one of the following holds.
{\rm (}Here let $\delta: S\to Y$ be the minimal resolution of $Y$.{\rm )}
\begin{itemize}
\item [\rm (2.0)]
$\kappa(S)=2$, $Y$ has at most canonical singularities with $h^{1}(\mathcal{O}_{Y})=0$ and $\chi(\mathcal{O}_{Y})=0$, and $H=K_{Y}+T$ with $H^{2}=1$, where $T$ is a non zero torsion divisor.
\item [\rm (2.1)]
$\kappa(S)=1$ and there exists an elliptic fibration $f:S\to C$ 
over a smooth curve $C$ such that $g(C)=1$, $\chi(\mathcal{O}_{S})=1$,
$q(S)=1$ and $\delta^{*}(H)F=1$, where $F$ is a general fiber of $f$.
In this case $Y$ has only rational singularities.
\item [\rm (2.2)]
$\kappa(S)=1$ and there exists an elliptic fibration $f:S\to C$ 
over a smooth curve $C$ such that $g(C)=0$, $\chi(\mathcal{O}_{S})=1$, $q(S)=0$
and one of the following holds.
{\rm (}Here let $t$ be the number of multiple fibers.{\rm )}
\begin{center}
\begin{tabular}{|c|c|c|c|}
$p_{g}(S)$ & $\delta^{*}(H)F$ & $t$ 
& $(m_{1}, \dots , m_{t})$ \\ \hline
$0$ & $6$ & $2$ & $(2,3)$ \\
$1$ & $4$ & $2$ & $(2,4)$ \\
$0$ & $3$ & $2$ & $(3,3)$ \\
$0$ & $2$ & $3$ & $(2,2,2)$ 
\end{tabular}
\end{center}
\item [\rm (2.3)]
$S$ is a one point blowing up of an Enriques surface $S^{\prime}$
and $\delta^{*}(H)=\mu^{*}(H^{\prime})-E_{\mu}$, where $\mu:S\to S^{\prime}$
is the blowing up at a point $P$, $H^{\prime}$ is an ample line bundle on $S^{\prime}$ and $E_{\mu}$ is the exceptional divisor.
\item [\rm (2.4)]
$\kappa(S)=-\infty$ and $q(S)=0$.
In this case $Y$ has only rational singularities.
\end{itemize}
\end{itemize}
\end{Theorem}
\noindent{\em Proof.}
Let $\delta: S\to Y$ be the minimal resolution of $Y$.
\\
(I) The case where $\chi(\mathcal{O}_{Y})\leq 0$.
\\
Then $h^{0}(m(K_{X}+L))\geq {m+1\choose 2}-(m-1)\chi(\mathcal{O}_{Y})$ by Theorem \ref{3-1T2}.
\\
Assume that $h^{0}(m(K_{X}+L))={m+1\choose 2}$.
Then, since $m\geq 2$, by the proof of Theorem \ref{3-1T2}, we have $\chi(\mathcal{O}_{Y})=0$, $H^{2}=1$ and $h^{0}(H)=1$.
Hence by Claim \ref{CL1}
\begin{eqnarray*}
1=h^{0}(H)&=&\chi(H)\\
&=&\chi(\mathcal{O}_{Y})+(1-g(Y,H))+H^{2}\\
&=&2-g(Y,H).
\end{eqnarray*}
Hence $g(Y,H)=1$.
Moreover since $\chi(\mathcal{O}_{Y})=0$, we have $h^{1}(\mathcal{O}_{Y})>0$.
Then $g(S,\delta^{*}(H))=g(Y,H)=1$.
In particular $\kappa(S)=-\infty$. Since $\delta^{*}(H)$ is nef and big, we have $g(S,\delta^{*}(H))\geq h^{1}(\mathcal{O}_{S})$ by \cite[Theorem 2.1]{Fukuma97-1}.
Moreover because $h^{1}(\mathcal{O}_{S})\geq h^{1}(\mathcal{O}_{Y})$, we have $1=g(S,\delta^{*}(H))\geq h^{1}(\mathcal{O}_{S})\geq h^{1}(\mathcal{O}_{Y})>0$.
Hence $g(S,\delta^{*}(H))=h^{1}(\mathcal{O}_{S})$ and $h^{1}(\mathcal{O}_{S})=h^{1}(\mathcal{O}_{Y})=1$.
Here we note that $\delta^{*}(H)$ is $\delta^{*}(H)$-minimal because $H$ is ample and $\delta$ is the minimal resolution.
Hence by \cite[Theorem 3.1]{Fukuma97-1}, we see that $(S,\delta^{*}(H))$ is a scroll over a smooth curve.
Then we can prove the following.
\begin{Claim}\label{CL2}
$\delta$ is the identity map.
\end{Claim}
\noindent
{\em Proof.}
Since $h^{1}(\mathcal{O}_{S})=h^{1}(\mathcal{O}_{Y})$, we see that $Y$ has the Albanese mapping by Lemma \ref{4L-T4}.
Then there exists an elliptic curve $C$ and morphisms $\alpha: Y\to C$ and $\alpha^{\prime}: S\to C$ such that $\alpha^{\prime}=\alpha\circ \delta$.
Here we note that $\alpha$ and $\alpha^{\prime}$ have connected fibers.
Since $\alpha^{\prime}$ is a $\mathbb{P}^{1}$-bundle over $C$, we see that any fiber of $\alpha^{\prime}$ is irreducible.
Assume that $\delta$ is not the identity map.
Then $\mbox{Sing}(Y)\neq\emptyset$ and $\alpha^{\prime}$ has non-irreducible fiber.
But this is a contradiction.
Therefore $\delta$ is the identity map. $\Box$
\\
\par
Hence $S\cong Y$, that is, $Y$ is smooth, and $(Y,H)$ is a scroll over a smooth elliptic curve $C$.
In particular, there exists an ample vector bundle $\mathcal{E}$ on $C$ such that $Y=\mathbb{P}_{C}(\mathcal{E})$ and $H=H(\mathcal{E})$.
Then $c_{1}(\mathcal{E})=1$ because $H^{2}=1$.
Therefore we see that $\mathcal{E}$ is an indecomposable ample vector bundle on $C$.
\\
\\
(II) Assume that $\chi(\mathcal{O}_{Y})>0$.
\\
Then we have $h^{0}(m(K_{X}+L))\geq {m\choose 2}+1$.
We consider $(X,L)$ with 
$h^{0}(m(K_{X}+L))={m\choose 2}+1$.
Then, since $m\geq 2$, by the proof of Theorem \ref{3-1T2} we obtain 
$\chi(\mathcal{O}_{Y})=\chi_{0}(Y,H)=1$, $\chi_{1}(Y,H)+\chi_{2}(Y,H)=0$ and $H^{2}=\chi_{2}(Y,H)=1$.
Hence we have $g(Y,H)=1-\chi_{1}(Y,H)=2$.
\par
Hence we see that a quasi-polarized surface $(S,\delta^{*}(H))$ 
is $\delta^{*}(H)$-minimal with $g(S,\delta^{*}(H))=2$
(Here we note that quasi-polarized surfaces of this type was studied in \cite{Bi-Fa-La06}.)
Here we note that $\delta^{*}(H)^{2}=1$ and $K_{S}\delta^{*}(H)=1$.
\par
Next we study $(S,\delta^{*}(H))$ with $g(S,\delta^{*}(H))=2$.
\\
(II.a) Assume that $\kappa(S)=2$.
Since $(\delta^{*}H)^{2}=H^{2}=1$ and $\delta^{*}(H)K_{S}=HK_{Y}=1$,
we see that $S$ is minimal because $(S,\delta^{*}(H))$ is $\delta^{*}(H)$-minimal (see Definition \ref{T-D1}).
By the Hodge index theorem we have $\delta^{*}(H)\equiv K_{S}$ 
and $K_{S}^{2}=1$.
Then $h^{1}(\mathcal{O}_{S})=0$ and $h^{1}(\mathcal{O}_{Y})=0$.
On the other hand $K_{S}=\delta^{*}(K_{Y})+E_{\delta}$ holds, where $E_{\delta}$ is a $\delta$-exceptional divisor.
Here we note that $E_{\delta}$ is not always effective.
Hence $\delta^{*}(H-K_{Y})\equiv E_{\delta}$.
If $E_{\delta}\neq 0$, then $(E_{\delta})^{2}<0$ by Grauert's criterion (e.g. \cite[(2.1) Theorem in Chapter III]{BaHuPeVa04}).
But since $\delta^{*}(H-K_{S})E_{\delta}=0$, this is impossible.
Therefore we have $E_{\delta}=0$ and $K_{S}=\delta^{*}(K_{Y})$.
Therefore  $Y$ has at most canonical singularities.
Namely the singularities of $Y$ are at most rational double points.
Therefore $Y$ is Gorenstein and $K_{Y}$ is a Cartier divisor.
Since $\delta^{*}(H)\equiv \delta^{*}(K_{Y})$, we have $H\equiv K_{Y}$.
If $H=K_{Y}$, then $h^{2}(H)=h^{2}(K_{Y})=h^{0}(\mathcal{O}_{Y})=1$.
But this contadicts Claim \ref{CL1}.
Therefore $H=K_{Y}+T$, where $T$ is a torsion divisor.
\\
(II.b) Next we consider the case where $\kappa(S)=1$.
Here we use the results of \cite{LaTu98}.
Let $h:S\to C$ be its elliptic fibration.
Then, since $(\delta^{*}H)^{2}=1$ and $K_{S}\delta^{*}H=1$, the following are possible from \cite{LaTu98}.
\begin{itemize}
\item [(1)] $h$ has no multiple fibers (see \cite[Table 3.1]{LaTu98}).
\begin{itemize}
\item [\rm (1.1)] $g(C)=0$, $\chi(\mathcal{O}_{S})=3$, $q(S)=0$, $p_{g}(S)=2$ and $\delta^{*}(H)F=1$. 
\item [\rm (1.2)] $g(C)=1$, $\chi(\mathcal{O}_{S})=1$, $q(S)=1$, $p_{g}(S)=1$ and $\delta^{*}(H)F=1$. (This is the type (2.1) in Theorem \ref{3-1T4}.)
\end{itemize}
\item [\rm (2)] The case where \cite[Table 4.1]{LaTu98}. (This is the type (2.2) in Theorem \ref{3-1T4}.)
\item [\rm (3)] $h$ has only one multiple fiber and its multiplicity is $2$.
In this case $g(C)=1$, $\chi(\mathcal{O}_{S})=0$, $q(S)=1$, $p_{g}(S)=0$ and $\delta^{*}HF=2$ (see the first case of \cite[Table 5.1]{LaTu98}).
\item [\rm (4)] The case where \cite[Table 5.2]{LaTu98}.
\end{itemize}

\begin{Lemma}\label{L3}
The cases {\rm (1.1)}, {\rm (3)} and {\rm (4)} above are impossible.
\end{Lemma}
\noindent{\em Proof.}
First we consider the case of (1.1).
In this case $\chi(\mathcal{O}_{S})=3>1=\chi(\mathcal{O}_{Y})$.
But this is impossible by Proposition \ref{2P-T1} because $\chi(\mathcal{O}_{Y})=\chi(\mathcal{O}_{X})$.
\par
Next we consider the case (3) above.
Since $q(S)=1$, $S$ has the Albanese fibration $\alpha: S\to B$, where $B$ is an elliptic curve.
In this case, since $C$ is also an elliptic curve, by the universality of the Albanese map we see that there exists a morphism $\lambda: B\to C$ such that $h=\lambda\circ \alpha$.
Because $\alpha$ and $h$ have connected fibers, we see that $\lambda$ is an isomorphism.
Namely we may assume that $\alpha=h$.
Moreover by Lemma \ref{4L-T4} the Albanese map of $Y$ can be defined, and let $\alpha_{Y}:Y\to B$ be its morphism.
But here $h$ is a quasi-bundle, so $\alpha$ is also a quasi-bundle.
 (For the definition of quasi-bundle, see \cite[Definition 1.1]{Serrano91}.)
Hence $\delta$ is an isomorphism because $\alpha=\alpha_{Y}\circ \delta$.
Therefore $Y\cong S$.
But then $\chi(\mathcal{O}_{Y})=\chi(\mathcal{O}_{S})=0$ and this is a contradiction.
\par
Finally we consider the case where (4).
Then by \cite[Proposition 5.1]{LaTu98}, $\delta^{*}H$ is ample.
Namely $\delta$ is an isomorphism.
But then $\chi(\mathcal{O}_{Y})=\chi(\mathcal{O}_{S})=0$ and
this is also impossible.
\par
This completes the proof of Lemma \ref{L3}. $\Box$
\\
(II.c) Next we consider the case where $\kappa(S)=0$.
Let $\mu:S\to S^{\prime}$ be the minimalization of $S$.
If $\delta$ is an isomorphism, then $\chi(\mathcal{O}_{S})=\chi(\mathcal{O}_{Y})=1$ and $S^{\prime}$ is an Enriques surface.
If $\delta$ is not an isomorphism, then 
since $g(S,\delta^{*}(H))=2$, by \cite[Proposition 3.2]{Bi-Fa-La06}
we see that $S^{\prime}$ is either an Enriques surface or a K3-surface.
If $S^{\prime}$ is birationally equivalent to a K3-surface, then $\chi(\mathcal{O}_{S^{\prime}})=2$.
But by Proposition \ref{2P-T1} this is impossible because $\chi(\mathcal{O}_{Y})=1$ in this case.
Therefore $S^{\prime}$ is birationally equivalent to an Enriques surface.
\\
(II.d) Next we consider the case where $\kappa(S)=-\infty$.
By Proposition \ref{2P-T1} we see that $\chi(\mathcal{O}_{S})\leq \chi(\mathcal{O}_{Y})=1$.
Since $g(S,\delta^{*}(H))=2$, we have $q(S)\leq 2$ by \cite[Theorem 2.1]{Fukuma97-1}.
By Lemma \ref{4L-T4}, we have $q(Y)=q(S)$ and if $q(Y)\geq 1$, then there exist the Albanese map of $Y$, $\alpha_{Y}: Y\to \mbox{Alb}(Y)$, and a morphism $\beta : \mbox{Alb}(S)\to \mbox{Alb}(Y)$ such that $\alpha_{Y}\circ\delta=\beta\circ\alpha_{S}$ holds, where $\alpha_{S}: S\to \mbox{Alb}(S)$ is the Albanese map of $S$.
Then $\alpha_{S}(S)$ and $\alpha_{Y}(Y)$ are smooth curves and $\alpha_{S}$ and $\alpha_{Y}$ have connected fibers (see \cite[Lemma 2.4.5]{BeSo-Book}).
Hence $\alpha_{S}(S)\cong \alpha_{Y}(Y)$.
\\
(i) If $q(S)=2$, then $g(S,\delta^{*}(H))=q(S)$ implies that $(S,\delta^{*}(H))$ is a scroll over a smooth curve by \cite[Theorem 3.1]{Fukuma97-1}.
Here we note that $\delta$ is an isomorphism because $S$ is a $\mathbb{P}^{1}$-bundle over $\alpha_{S}(S)$.
But then $\chi(\mathcal{O}_{Y})=\chi(\mathcal{O}_{S})=-1$ and this is impossible.\\
(ii) Next we consider the case where $q(S)=1$.
Assume that $K_{S}+\delta^{*}(H)$ is not nef.
Then there exists an extremal rational curve $E$ on $S$ such that $(K_{S}+\delta^{*}(H))E<0$.
If $E$ is a $(-1)$-curve, then $(K_{S}+\delta^{*}(H))E\geq 0$
since $(S,\delta^{*}(H))$ is $\delta^{*}(H)$-minimal.
Hence $S$ is a $\mathbb{P}^{1}$-bundle over a smooth elliptic curve $C$ and $E$ is a fiber of this because $q(S)=1$.
Let $f: S\to C$ be its morphism.
Moreover we see that $\delta^{*}(H)F=1$ for any fiber $F$ of $f$ because $(K_{S}+\delta^{*}(H))F<0$.
Then $g(S,\delta^{*}(H))=q(S)=1$. 
But this contradicts to $g(S,\delta^{*}(H))=g(Y,H)=2$.
Hence $K_{S}+\delta^{*}(H)$ is nef.
So we get $0\leq (K_{S}+\delta^{*}(H))^{2}=K_{S}^{2}+2K_{S}\delta^{*}(H)+(\delta^{*}(H))^{2}=3+K_{S}^{2}$, that is, $-3\leq K_{S}^{2}$.
On the other hand $K_{S}^{2}\leq 0$ and $K_{S}^{2}=0$ if and only if $S$ is minimal.
Hence $S$ is at most three points blowing up of a $\mathbb{P}^{1}$-bundle over $C$.
\\
(ii.1) Assume that $S$ is a $\mathbb{P}^{1}$-bundle over $C$.
Then $S\cong Y$ because every exceptional curve of $\delta$ is contained in a fiber of $\alpha_{S}$.
But this is impossible because $\chi(\mathcal{O}_{S})=0\neq1=\chi(\mathcal{O}_{Y})$.
\\
(ii.2) Assume that $S$ is one point blowing up of a $\mathbb{P}^{1}$-bundle over $C$.
Then $S$ has one singular fiber $F_{1}$ and $F_{1}=C_{1}+C_{2}$, where each $C_{i}$ is a $(-1)$-curve and $C_{1}C_{2}=1$.
Since $\delta$ is the minimal resolution, we have $S\cong Y$.
But this is also impossible by the same reason as in (ii.1).
\\
(ii.3) Assume that $S$ is two point blowing up of a $\mathbb{P}^{1}$-bundle over $C$.
Then the following two cases possibly occur:\\
(ii.3.1) $\alpha_{S}$ has one singular fiber $F$ and $F=C_{1}+C_{2}+C_{3}$, where $C_{1}$ and $C_{3}$ are $(-1)$-curves and $C_{2}$ is a $(-2)$-curve such that $C_{1}C_{2}=1$, $C_{2}C_{3}=1$ and $C_{1}C_{3}=0$.
\\
(ii.3.2) $f$ has two singular fibers $F_{1}$ and $F_{2}$ such that $F_{1}=C_{1}+C_{2}$, $F_{2}=C_{3}+C_{4}$, where each $C_{i}$ is a $(-1)$-curve with $C_{1}C_{2}=1$ and $C_{3}C_{4}=1$.\\
By the same argument as (ii.2), (ii.3.2) cannot occur.
So we consider the case where (ii.3.1).
Then since $\delta$ is the minimal resolution, the exceptional curve of $\delta$ is $C_{2}$. So $Y$ is rational by Artin's criterion \cite[(3.2) Theorem in ChapterIII]{BaHuPeVa04}.
But this is impossible because $\chi(\mathcal{O}_{S})=0\neq1=\chi(\mathcal{O}_{Y})$.
\\
(ii.4) Assume that $S$ is three point blowing up of a $\mathbb{P}^{1}$-bundle over $C$.
Then the following four cases possibly occur:\\
(ii.4.1) $\alpha_{S}$ has one singular fiber $F$ and $F=C_{1}+C_{2}+C_{3}+C_{4}$, where $C_{2}$, $C_{3}$ and $C_{4}$ are $(-1)$-curves and $C_{1}$ is a $(-3)$-curve such that $C_{1}C_{i}=1$ for every $i$ with $i=2, 3, 4$, $C_{j}C_{k}=0$ with $j,k\in\{ 2, 3, 4\}$ and $j\neq k$.
\\
(ii.4.2) $\alpha_{S}$ has one singular fiber $F$ and $F=C_{1}+C_{2}+C_{3}+C_{4}$, where $C_{1}$ and $C_{4}$ are $(-1)$-curves, and $C_{2}$ and $C_{3}$ are $(-2)$-curves such that $C_{i}C_{i+1}=1$ for every $i$ with $i=1, 2, 3$, $C_{j}C_{k}=0$ with $|j-k|\geq 2$.
\\
(ii.4.3) $\alpha_{S}$ has two singular fibers $F_{1}$ and $F_{2}$ such that $F_{1}=C_{1}+C_{2}+C_{3}$, $F_{2}=C_{4}+C_{5}$, where $C_{i}$ is a $(-1)$-curve for every $i\neq 2$ and $C_{2}$ is a $(-2)$-curve such that $C_{1}C_{2}=1$, $C_{2}C_{3}=1$, $C_{1}C_{3}=0$ and $C_{4}C_{5}=1$.
\\
(ii.4.4) $f$ has three singular fibers $F_{1}$, $F_{2}$ and $F_{3}$ such that $F_{1}=C_{1}+C_{2}$, $F_{2}=C_{3}+C_{4}$ and $F_{3}=C_{5}+C_{6}$, where each $C_{i}$ is a $(-1)$-curve such that $C_{i}C_{i+1}=1$ with $i\in \{ 1, 3, 5\}$.
\\
By the same argument as above, 
in these 4 cases we see that $\delta$ is an isomorphism
or $Y$ has rational singularities.
But this is impossible because $\chi(\mathcal{O}_{S})=0\neq \chi(\mathcal{O}_{Y})$.
\par
Therefore the case where $q(S)=1$ cannot occur.
By the above argument, we see that $q(S)=0$.
Then $\chi(\mathcal{O}_{S})=1=\chi(\mathcal{O}_{Y})$ and by Proposition \ref{2P-T1} we have $h^{0}(R^{1}\delta_{*}(\mathcal{O}_{S}))=0$.
So $Y$ has rational singularities.
This completes the proof. $\Box$

\section{Main results}\label{S5}

Let $(X,L)$ be a polarized manifold of dimension $3$.
In this section, we consider $h^{0}(m(K_{X}+L))$.
First by Theorems \ref{3-1T1} and \ref{3-1T2} we have the following.

\begin{Theorem}\label{3-2T1}
Let $(X,L)$ be a polarized manifold of dimension $3$.
\begin{itemize}
\item [\rm (1)]
Assume that $\kappa(K_{X}+L)=0$.
Then $h^{0}(m(K_{X}+L))=1$ for every positive integer $m$.
\item [\rm (2)] Assume that $\kappa(K_{X}+L)=1$.
Then for every positive integer $m$ the following holds.
\[
h^{0}(m(K_{X}+L))\geq
\left\{
\begin{array}{lc}
(m-1)(h^{1}(\mathcal{O}_{X})-1)+mh^{1}(\mathcal{O}_{X}) & \mbox{if $h^{1}(\mathcal{O}_{X})\geq 1$,} \\
m+1 & \mbox{if $h^{1}(\mathcal{O}_{X})=0$.}
\end{array} \right. 
\]
\item [\rm (3)] Assume that $\kappa(K_{X}+L)=2$.
Then for every positive integer $m$ the following holds.
\[
h^{0}(m(K_{X}+L))\geq
\left\{
\begin{array}{lc}
{m+1\choose 2}-(m-1)\chi(\mathcal{O}_{X}) & \mbox{if $\chi(\mathcal{O}_{X})\leq 0$,} \\
{m\choose 2}+\chi(\mathcal{O}_{X}) & \mbox{if $\chi(\mathcal{O}_{X})>0$.}
\end{array} \right. 
\]
\end{itemize}
\end{Theorem}
\noindent{\em Proof.}
Let $(M,A)$ be a reduction of $(X,L)$.
Here we note that $h^{0}(m(K_{X}+L))=h^{0}(m(K_{M}+A))$ for any positive integer $m$.
If $\kappa(K_{X}+L)=0$, then $(M,A)$ is a Mukai manifold, that is, $\mathcal{O}_{M}(K_{M}+A)=\mathcal{O}_{M}$ by \cite[Theorem 7.5.3]{BeSo-Book}.
This implies that $h^{0}(m(K_{X}+L))=h^{0}(m(K_{M}+A))=1$.
\par
If $\kappa(K_{X}+L)=1$ (resp. $2$), then by \cite[Theorem 7.5.3]{BeSo-Book} there exist a smooth projective curve $C$ (resp. a normal projective surface $Y$), and a fiber space $f:M\to C$ (resp. $M\to Y$) such that $K_{M}+A=f^{*}(H)$ for some ample line bundle $H$ on $C$ (resp. $Y$).
Moreover we have $h^{1}(\mathcal{O}_{X})=h^{1}(\mathcal{O}_{M})=h^{1}(\mathcal{O}_{C})$ (resp. $h^{i}(\mathcal{O}_{X})=h^{i}(\mathcal{O}_{M})=h^{i}(\mathcal{O}_{Y})$ for $i=0, 1, 2$ and $h^{3}(\mathcal{O}_{X})=0$).
Hence by Theorems \ref{3-1T1} and \ref{3-1T2} we get the assertion. $\Box$
\\
\par
Next we consider the case where $\kappa(K_{X}+L)=3$.
Then the following is obtained.

\begin{Theorem}\label{3-2T2}
Let $(X,L)$ be a polarized manifold of dimension $3$.
Assume that $\kappa(K_{X}+L)=3$.
Then we have
\[
h^{0}(m(K_{X}+L))
\geq \left\{
\begin{array}{lc}
\frac{1}{8}m^{3}+\frac{1}{4}m^{2}+1 & \mbox{if $m$ is even with $m\geq 2$,} \\
\frac{1}{8}m^{3}+\frac{1}{4}m^{2}+\frac{1}{8}m+1 & \mbox{if $m$ is odd with $m\geq 3$.}
\end{array} \right.
\]
\end{Theorem}
\noindent{\em Proof.}
Let $(M,A)$ be a reduction of $(X,L)$.
By assumption and \cite[Proposition 7.6.9]{BeSo-Book} we see that $K_{M}+A$ is nef.
\\
(I) The case where $m$ is even with $m\geq 2$.
\\
Then by Proposition \ref{SKB} we have the following.
\begin{eqnarray*}
h^{0}(m(K_{X}+L))
&=&h^{0}(m(K_{M}+A))\\
&=&h^{0}\left(\left(\frac{m}{2}+1\right)K_{M}+\frac{m}{2}A\right)
+g_{2}\left(M,\left(\frac{m}{2}-1\right)(K_{M}+A)+A\right)\\
&&\ \ \ -h^{1}(\mathcal{O}_{M})
+g_{1}\left(M,\left(\frac{m}{2}-1\right)(K_{M}+A)+A, \frac{m}{2}(K_{M}+A)\right).
\end{eqnarray*}

Since $((m/2)-1)(K_{M}+A)+A$ is ample and $\kappa(K_{M}+((m/2)-1)(K_{M}+A)+A)=\kappa(K_{M}+A)=3$, we have 
$g_{2}(M,((m/2)-1)(K_{M}+A)+A)\geq h^{1}(\mathcal{O}_{M})$ by Theorem \ref{Theorem1.4}.
On the other hand, by Remark \ref{B6} (2) we have
\begin{eqnarray*}
&&g_{1}\left(M,\left(\frac{m}{2}-1\right)(K_{M}+A)+A, \frac{m}{2}(K_{M}+A)\right)\\
&&=1+\frac{1}{2}\left(K_{M}+\left(\frac{m}{2}-1\right)(K_{M}+A)+A+\frac{m}{2}(K_{M}+A)\right)\\
&&\ \ \ \times\left(\left(\frac{m}{2}-1\right)(K_{M}+A)+A)\right)\left(\frac{m}{2}(K_{M}+A)\right)\\
&&=1+\frac{m^{2}(m-2)}{8}(K_{M}+A)^{3}+\frac{m^{2}}{4}(K_{M}+A)^{2}A.
\end{eqnarray*}

We also note that $(K_{M}+A)^{3}\geq 1$ and $(K_{M}+A)^{2}A\geq 1$.
\\
If $(K_{M}+A)^{2}A=1$, then by Proposition \ref{GHIT} we see that
$(K_{M}+A)A^{2}=1$ and $A^{3}=1$ because $(K_{M}+A)A^{2}>0$.
Hence $g_{1}(M,A)=2$.
Therefore by \cite[(1.10) Theorem and Section 2]{Fujita87-2} we see that $K_{M}=\mathcal{O}_{M}$ and $h^{0}(A)\geq 1$ since $\kappa(K_{M}+A)=3$.
On the other hand, we have
$$h^{0}(m(K_{M}+A))=h^{0}(mA)=\chi(mA)=\frac{1}{6}m^{3}A^{3}+\frac{1}{12}mc_{2}(M)A$$
because $h^{i}(mA)=h^{i}(K_{M}+mA)=0$ for every $i>0$.
Since $h^{0}(A)\geq 1$, we get 
$$1\leq h^{0}(A)=\frac{1}{6}A^{3}+\frac{1}{12}c_{2}(M)A.$$
Hence $(1/12)c_{2}(M)A\geq 1-(1/6)A^{3}=5/6$.
So we obtain
\begin{eqnarray*}
h^{0}(m(K_{M}+A))
&=&\frac{1}{6}m^{3}A^{3}+\frac{1}{12}mc_{2}(M)A\\
&\geq& \frac{1}{6}m^{3}+\frac{5}{6}m.
\end{eqnarray*}
\noindent
\\
If $(K_{M}+A)^{2}A\geq 2$, then
\begin{eqnarray*}
h^{0}(m(K_{M}+A))
&\geq & 1+\frac{m^{2}(m-2)}{8}+2\frac{m^{2}}{4}\\
&=&\frac{1}{8}m^{3}+\frac{1}{4}m^{2}+1.
\end{eqnarray*}

Here we note that
$(1/6)m^{3}+(5/6)m-((1/8)m^{3}+(1/4)m^{2}+1)=(1/24)(m-2)((m-2)^{2}+8)\geq 0$.
So if $m$ is even with $m\geq 2$, then we have $h^{0}(m(K_{M}+A))\geq (1/8)m^{3}+(1/4)m^{2}+1$.
\\
(II) The case where $m$ is odd with $m\geq 3$.
\\
Here we use the following equality which is obtained from Proposition \ref{SKB}.
\begin{eqnarray*}
h^{0}(m(K_{X}+L))
&=&h^{0}(m(K_{M}+A))\\
&=&h^{0}\left(\left(\frac{m+1}{2}+1\right)K_{M}+\frac{m+1}{2}A\right)\\
&&\ \ \ +g_{2}\left(M,\left(\frac{m-1}{2}-1\right)(K_{M}+A)+A\right)-h^{1}(\mathcal{O}_{M})\\
&&\ \ \ +g_{1}\left(M,\left(\frac{m-1}{2}-1\right)(K_{M}+A)+A, \frac{m+1}{2}(K_{M}+A)\right).
\end{eqnarray*}
Since $(-1+(m-1)/2)(K_{M}+A)+A$ is ample and $\kappa(K_{M}+(-1+(m-1)/2)(K_{M}+A)+A)=\kappa(((m-1)/2)(K_{M}+A))=\kappa(K_{M}+A)=3$, we have 
$g_{2}(M,(-1+(m-1)/2)(K_{M}+A)+A)\geq h^{1}(\mathcal{O}_{M})$ by Theorem \ref{Theorem1.4}.
On the other hand, 
\begin{eqnarray*}
&&g_{1}\left(M,\left(\frac{m-1}{2}-1\right)(K_{M}+A)+A, \frac{m+1}{2}(K_{M}+A)\right)\\
&&=1+\frac{1}{2}\left(K_{M}+\left(\frac{m-1}{2}-1\right)(K_{M}+A)+A+\frac{m+1}{2}(K_{M}+A)\right)\\
&&\ \ \ \times\left(\left(\frac{m-1}{2}-1\right)(K_{M}+A)+A\right)\left(\frac{m+1}{2}(K_{M}+A)\right)\\
&&=1+\frac{m(m+1)(m-3)}{8}(K_{M}+A)^{3}+\frac{m(m+1)}{4}(K_{M}+A)^{2}A.
\end{eqnarray*}
If $(K_{M}+A)^{2}A=1$, then by the same argument as above we see that
$$h^{0}(m(K_{M}+A))\geq\frac{1}{6}m^{3}+\frac{5}{6}m.$$
\noindent
If $(K_{M}+A)^{2}A\geq 2$, then we have
\begin{eqnarray*}
h^{0}(m(K_{M}+A))
&\geq & 1+\frac{m(m+1)(m-3)}{8}+\frac{m(m+1)}{2}\\
&=&\frac{1}{8}m^{3}+\frac{1}{4}m^{2}+\frac{1}{8}m+1.
\end{eqnarray*}
Here we note that
$(1/6)m^{3}+(5/6)m-((1/8)m^{3}+(1/4)m^{2}+(1/8)m+1)=(1/24)(m-3)((m-(3/2))^{2}+23/4)\geq 0$.
So if $m$ is odd with $m\geq 3$, then we have $h^{0}(m(K_{M}+A))\geq (1/8)m^{3}+(1/4)m^{2}+(1/8)m+1$.
This completes the proof of Theorem \ref{3-2T2}. $\Box$

\begin{Remark}\label{R1}
By Theorem \ref{3-2T2} we see that if $\kappa(K_{X}+L)=3$, then for every integer $m$ with $m\geq 2$, we have
$$h^{0}(m(K_{X}+L))\geq \frac{1}{8}m^{3}+\frac{1}{4}m^{2}+1.$$
\end{Remark}

If $\kappa(K_{X}+L)=3$ and $m=2$, then by Theorem \ref{3-2T2} or \cite[Theorem 5.4 (2)]{Fukuma08-3} we have $h^{0}(2(K_{X}+L))\geq 3$.
So it is interesting to study $(X,L)$ with $\kappa(K_{X}+L)=3$ and small $h^{0}(2(K_{X}+L))$.
The following results (Theorems \ref{EC1} and \ref{EC2}) give a classification of these $(X,L)$.
\\
\par
First we note the following which will be used later.

\begin{Proposition}\label{32P-T1}
Let $(X,L)$ be a polarized manifold of dimension $3$.
Then the following equalities holds.
\begin{eqnarray}
&&h^{0}(2K_{X}+2L)-h^{0}(2K_{X}+L) \label{m2.3.1}\\
&&=g_{2}(X,L)-h^{1}(\mathcal{O}_{X})+g_{1}(X,K_{X}+L,L), \nonumber \\
&&h^{0}(2K_{X}+2L)-h^{0}(K_{X}+L) \label{m2.3.2}\\
&&=g_{2}(X,K_{X}+L)-h^{1}(\mathcal{O}_{X})+g_{1}(X,K_{X}+L,L). \nonumber 
\end{eqnarray}
\end{Proposition}
\noindent
{\em Proof.} These equalities are obtained from Proposition \ref{SKB}. $\Box$

\begin{Notation}\label{32N-T1}
Let $(X,L)$ be a polarized manifold of dimension $3$ and let $(M,A)$ be a reduction of $(X,L)$.
Set $d_{1}:=g_{2}(M,A)-h^{1}(\mathcal{O}_{M})$ and $d_{2}:=g_{2}(M,K_{M}+A)-h^{1}(\mathcal{O}_{M})$.
Then we see that
\begin{eqnarray*}
d_{2}-d_{1}
&=&\frac{1}{12}(K_{M}+A)(6K_{M}+6A)K_{M}+\frac{1}{12}c_{2}(M)K_{M}\\
&=&\frac{1}{12}(K_{M}+A)(6K_{M}+6A)K_{M}-2\chi(\mathcal{O}_{M}).
\end{eqnarray*}
Therefore
\begin{equation}
d_{2}-d_{1}+2\chi(\mathcal{O}_{M})=\frac{1}{2}(K_{M}+A)^{2}K_{M}.\label{m2.3.7}
\end{equation}
\end{Notation}

\begin{Theorem}\label{EC1}
Let $(X,L)$ be a polarized manifold of dimension $3$.
Assume that $\kappa(K_{X}+L)=3$.
Then $h^{0}(2(K_{X}+L))=3$ if and only if $(X,L)$ satisifes $L^{3}=1$, 
$\mathcal{O}_{X}(K_{X})=\mathcal{O}_{X}$, 
$h^{1}(\mathcal{O}_{X})=0$ and $h^{0}(L)=1$.
\end{Theorem}
\noindent{\em Proof.}
($\alpha$) Assume that $h^{0}(2(K_{X}+L))=3$.
\\
Let $(M,A)$ be a reduction of $(X,L)$.
Then by assumption we see that $K_{M}+A$ is nef and big.
First we prove the following claim.

\begin{Claim}\label{CL4}
$h^{0}(K_{M}+A)\leq 2$.
\end{Claim}
\noindent
{\em Proof.}
Assume that $h^{0}(K_{M}+A)\geq 3$.
Then by Lemma \ref{Lemma B} 
we have $h^{0}(2(K_{M}+A))\geq 2h^{0}(K_{M}+A)-1\geq 5$.
This is a contradiction. $\Box$
\\
\par
By Proposition \ref{32P-T1} (\ref{m2.3.1}) and Theorem \ref{Theorem1.4}, we see that
\begin{eqnarray*}
3=h^{0}(2K_{M}+2A)
&\geq & g_{2}(M,A)-h^{1}(\mathcal{O}_{M})+g_{1}(M,K_{M}+A,A)\\
&\geq & g_{1}(M,K_{M}+A,A)\\
&=& 1+(K_{M}+A)^{2}A.
\end{eqnarray*}

Hence we have $(K_{M}+A)^{2}A\leq 2$.
On the other hand, since $1\leq (K_{M}+A)^{2}A$ we get

\begin{equation}
1\leq (K_{M}+A)^{2}A\leq 2. \label{m2.3.3}
\end{equation}
Namely the following holds.

\begin{equation}
2\leq g_{1}(M,K_{M}+A,A)\leq 3. \label{m2.3.4}
\end{equation}

Since $g_{1}(M,K_{M}+A,A)\leq 3$, by Proposition \ref{32P-T1} (\ref{m2.3.2})
we get
$h^{0}(2(K_{M}+A))-h^{0}(K_{M}+A)\leq g_{2}(M, K_{M}+A)-h^{1}(\mathcal{O}_{M})+3$.
By Claim \ref{CL4} and $h^{0}(2(K_{M}+A))=3$, we see that 
\begin{eqnarray*}
3-2
&\leq&h^{0}(2(K_{M}+A))-h^{0}(K_{M}+A)\\
&\leq&g_{2}(M,K_{M}+A)-h^{1}(\mathcal{O}_{M})+3.
\end{eqnarray*}
Namely,
\begin{eqnarray}
g_{2}(M, K_{M}+A)-h^{1}(\mathcal{O}_{M})&\geq& -2. \label{m2.3.5}
\end{eqnarray}

From Proposition \ref{32P-T1} (\ref{m2.3.2}), (\ref{m2.3.4}) and the assumption that $h^{0}(2(K_{M}+A))=3$, we have
\begin{eqnarray*}
3&\geq&h^{0}(2(K_{M}+A))-h^{0}(K_{M}+A)\\
&=&g_{2}(M,K_{M}+A)-h^{1}(\mathcal{O}_{M})+g_{1}(M,K_{M}+A,A) \\
&\geq&g_{2}(M,K_{M}+A)-h^{1}(\mathcal{O}_{M})+2.
\end{eqnarray*}
Hence we have
\begin{eqnarray}
g_{2}(M, K_{M}+A)-h^{1}(\mathcal{O}_{M})&\leq& 1. \label{m2.3.6}
\end{eqnarray}

By (\ref{m2.3.4}) and Proposition \ref{32P-T1} (\ref{m2.3.1}), we have 
\begin{eqnarray*}
3&\geq&h^{0}(2(K_{M}+A))-h^{0}(2K_{M}+A)\\
&=&g_{2}(M,A)-h^{1}(\mathcal{O}_{M})+g_{1}(M,K_{M}+A,A) \\
&\geq&g_{2}(M,A)-h^{1}(\mathcal{O}_{M})+2.
\end{eqnarray*}
Hence $1\geq g_{2}(M,A)-h^{1}(\mathcal{O}_{M})$.
From this and Theorem \ref{Theorem1.4} we have
\begin{equation}
d_{1}=0, 1. \label{m2.3.9}
\end{equation}
We also note that 
\begin{equation}
-2\leq d_{2}\leq 1 \label{m2.3.10}
\end{equation} 
by (\ref{m2.3.5}) and (\ref{m2.3.6}).
\\
\\
(I) If $(K_{M}+A)^{2}A=1$, then $(K_{M}+A)A^{2}=1$ and $A^{3}=1$ by Proposition \ref{GHIT}.
Therefore we get $g_{1}(M,A)=2$.
Since $\kappa(K_{M}+A)=3$, by  \cite[(1.10) Theorem and Section 2]{Fujita87-2} we see that $K_{M}=\mathcal{O}_{M}$, $h^{1}(\mathcal{O}_{M})=0$ and $h^{0}(A)=1$.
By the Riemann-Roch theorem we get $\chi(tA)=(1/6)A^{3}t^{3}+(1/12)c_{2}(M)At$.
Since $h^{0}(2K_{M}+2A)=\chi(2K_{M}+2A)=\chi(2A)$,
we get $h^{0}(2K_{M}+2A)=(4/3)A^{3}+(1/6)c_{2}(M)A$.
Therefore $3=h^{0}(2K_{M}+2A)=(4/3)A^{3}+(1/6)c_{2}(M)A=(4/3)+(1/6)c_{2}(M)A$.
Namely $c_{2}(M)A=10$.
Here we note that $(M,A)\cong (X,L)$ because $A^{3}=1$.
\\
\\
(II) Next we assume that 
\begin{equation}
(K_{M}+A)^{2}A=2. \label{m2.3.11}
\end{equation}
We will prove that this case cannot occur.
Since $(K_{M}+A)^{2}A=2$, by Proposition \ref{GHIT} we have 
\begin{equation}
1\leq (K_{M}+A)^{3}\leq 4. \label{m2.3.12}
\end{equation}
By using (\ref{m2.3.7}), (\ref{m2.3.9}), (\ref{m2.3.10}), (\ref{m2.3.11}) and (\ref{m2.3.12}), we can determine the value of $\chi(\mathcal{O}_{M})$. 
For example, assume that $d_{1}=0$ and $d_{2}=-2$.
Then $d_{2}-d_{1}=-2$ and $(K_{M}+A)^{2}K_{M}=4\chi(\mathcal{O}_{M})-4$ by (\ref{m2.3.7}).
Since $(K_{M}+A)^{2}A=2$, we have $(K_{M}+A)^{3}=4\chi(\mathcal{O}_{M})-2$.
By considering (\ref{m2.3.12}) we have $\chi(\mathcal{O}_{M})=1$.
By the same argument as this, 
we can get the following list:

\begin{center}
\begin{tabular}{cccccc}\hline
$d_{1}$ & $d_{2}$ & $d_{2}-d_{1}$ & $(K_{M}+A)^{2}K_{M}$ & $(K_{M}+A)^{3}$ & $\chi(\mathcal{O}_{M})$\\
\hline
$0$ & $-2$ & $-2$ & $4\chi(\mathcal{O}_{M})-4$ & $4\chi(\mathcal{O}_{M})-2$ & 
$1$ \\
$0$ & $-1$ & $-1$ & $4\chi(\mathcal{O}_{M})-2$ & $4\chi(\mathcal{O}_{M})$ & $1$ \\
$0$ & $0$ & $0$ & $4\chi(\mathcal{O}_{M})$ & $4\chi(\mathcal{O}_{M})+2$ & $0$ \\
$0$ & $1$ & $1$ & $4\chi(\mathcal{O}_{M})+2$ & $4\chi(\mathcal{O}_{M})+4$ & $0$ \\
$1$ & $-2$ & $-3$ & $4\chi(\mathcal{O}_{M})-6$ & $4\chi(\mathcal{O}_{M})-4$ & 
$2$ \\
$1$ & $-1$ & $-2$ & $4\chi(\mathcal{O}_{M})-4$ & $4\chi(\mathcal{O}_{M})-2$ & $1$ \\
$1$ & $0$ & $-1$ & $4\chi(\mathcal{O}_{M})-2$ & $4\chi(\mathcal{O}_{M})$ & $1$ \\
$1$ & $1$ & $0$ & $4\chi(\mathcal{O}_{M})$ & $4\chi(\mathcal{O}_{M})+2$ & $0$
\end{tabular}
\end{center}

By this list, we see that $(K_{M}+A)^{3}=2$ or $4$.
\par
Assume that $(K_{M}+A)^{3}=4$.
Then by Proposition \ref{GHIT} we have 
\begin{eqnarray*}
4&=&((K_{M}+A)^{2}A)^{2} \\
&\geq & ((K_{M}+A)^{3})((K_{M}+A)A^{2}) \\
&\geq &4(K_{M}+A)A^{2}.
\end{eqnarray*}
Since $K_{M}+A$ is nef and big, we see that $(K_{M}+A)A^{2}\geq 1$.
Therefore $(K_{M}+A)A^{2}=1$.
But by Proposition \ref{GHIT}, we have 
$1=((K_{M}+A)A^{2})^{2}\geq ((K_{M}+A)^{2}A)A^{3}=2A^{3}\geq 2$,
and this is impossible.
\par
Assume that $(K_{M}+A)^{3}=2$.
Then by Proposition \ref{GHIT} we have 
\begin{eqnarray*}
4&=&((K_{M}+A)^{2}A)^{2} \\
&\geq & ((K_{M}+A)^{3})((K_{M}+A)A^{2}) \\
&=&2(K_{M}+A)A^{2}.
\end{eqnarray*}
Hence we have 
$(K_{M}+A)A^{2}\leq 2$.
By Proposition \ref{GHIT} we see that 
$((K_{M}+A)A^{2})^{2}\geq ((K_{M}+A)^{2}A)A^{3}=2A^{3}\geq 2$.
Therefore $(K_{M}+A)A^{2}=2$ and $A^{3}\leq 2$ because $(K_{M}+A)A^{2}\geq 1$.
But since $(K_{M}+A)A^{2}=2$, we have $A^{3}=2$ because $(K_{M}+2A)A^{2}$ is even.
Therefore 
$((K_{M}+A)A^{2})^{2}=4=((K_{M}+A)^{2}A)A^{3}$ holds and
$K_{M}+A\equiv A$ by \cite[Corollary 2.5.4]{BeSo-Book} since $A$ is ample.
Namely $K_{M}\equiv 0$.
Now since $g_{1}(M,A,K_{M}+A)=1+(K_{M}+A)^{2}A=3$, we see that $h^{0}(2K_{M}+A)=-d_{1}$ by Proposition \ref{32P-T1} (\ref{m2.3.1}).
Since $d_{1}=0$ or $1$ by (\ref{m2.3.9}), we have $d_{1}=0$.
On the other hand, $h^{i}(K_{M}+K_{M}+A)=0$ for every integer $i$ with $i>0$ because $K_{M}+A$ is nef and big.
So by the Riemann-Roch theorem we have $h^{0}(2K_{M}+A)=\chi(2K_{M}+A)=\chi(A)=(1/6)A^{3}+(1/12)c_{2}(M)A$.
Since $A^{3}=2$, we have $c_{2}(M)A=-4$ if $d_{1}=0$.
Here we calculate $h^{0}(2(K_{M}+A))$.
Since $K_{M}+2A$ is ample, then $h^{i}(2K_{M}+2A)=0$ for $i>0$.
Therefore
\begin{eqnarray*}
h^{0}(2(K_{M}+A))
&=&\chi(2(K_{M}+A))\\
&=&\chi(2A)\\
&=&\frac{4}{3}A^{3}+\frac{1}{6}c_{2}(M)A\\
&=&2.
\end{eqnarray*}
But this is impossinble because we assume that $h^{0}(2(K_{M}+A))=3$.
\\
($\beta$) Assume that $(X,L)$ satisfies $L^{3}=1$, 
$\mathcal{O}_{X}(K_{X})=\mathcal{O}_{X}$, 
$h^{1}(\mathcal{O}_{X})=0$ and $h^{0}(L)=1$.
Then $h^{0}(2K_{X}+L)=h^{0}(K_{X}+L)=h^{0}(L)=1$ and $h^{2}(\mathcal{O}_{X})=h^{1}(K_{X})=h^{1}(\mathcal{O}_{X})=0$.
Hence $g_{2}(X,L)=h^{0}(K_{X}+L)-h^{0}(K_{X})+h^{2}(\mathcal{O}_{X})=0$.
Moreover $g_{1}(X,K_{X}+L,L)=1+L^{3}=2$.
Therefore by Proposition \ref{32P-T1} (\ref{m2.3.1}) we have
\begin{eqnarray*}
h^{0}(2(K_{X}+L))
&=&h^{0}(2K_{X}+L)+g_{2}(X,L)-h^{1}(\mathcal{O}_{X})+g_{1}(X,K_{X}+L,L)\\
&=&3.
\end{eqnarray*}
This completes the proof. $\Box$

\begin{Remark}\label{}
(i) By Theorem \ref{EC1}, we see that if $\kappa(K_{X}+L)=3$ and $h^{0}(2(K_{X}+L))=3$, then $h^{0}(K_{X}+L)=1$.\\
(ii) There exists an example of $(X,L)$ which satisfies $\kappa(K_{X}+L)=3$ and $h^{0}(2K_{X}+2L)=3$.
See \cite[Example 3.1 (4)]{Fukuma10}.
\end{Remark}

Next we consider the case where $(X,L)$ satisfies $\kappa(K_{X}+L)=3$ and $h^{0}(2K_{X}+2L)=4$.

\begin{Theorem}\label{EC2}
Let $(X,L)$ be a polarized manifold of dimension $3$ 
and let $(M,A)$ be a reduction of $(X,L)$.
Assume that $\kappa(K_{X}+L)=3$.
Then $h^{0}(2(K_{X}+L))=4$ if and only if $(M,A)$ is one of the following.
\begin{itemize}
\item [\rm (1)] $K_{M}\equiv 0$, $A^{3}=2$, $\chi(\mathcal{O}_{M})=0$ and $h^{0}(A)=1$.
\item [\rm (2)] $(K_{M}+A)^{2}A=3$, $(K_{M}+A)^{3}=1$, $g_{2}(M,A)=h^{1}(\mathcal{O}_{M})=1$, $h^{2}(\mathcal{O}_{M})=0$ , $h^{3}(\mathcal{O}_{M})=0$ and $(M,K_{M}+A)$ is birationally equivalent to a scroll over an elliptic curve.
\end{itemize}
\end{Theorem}
\noindent{\em Proof.}
($\alpha$) Assume that $h^{0}(2(K_{X}+L))=4$.
\\
First we prove the following claim.
\begin{Claim}\label{T-CL1}
One of the following holds:
\begin{itemize}
\item [\rm (i)] $g(M,A)=2$.
\item [\rm (ii)] $(M,A)$ satisfies {\rm (1)} in Theorem {\rm \ref{EC2}}.
\item [\rm (ii)] $(M,A)$ satisfies {\rm (2)} in Theorem {\rm \ref{EC2}}.
\end{itemize}
\end{Claim}
\noindent
{\em Proof.}
If $h^{0}(K_{M}+A)\geq 3$, then 
by Lemma \ref{Lemma B}
we see that $h^{0}(2K_{M}+2A)\geq 2h^{0}(K_{M}+A)-1\geq 5$ and this is impossible.
Hence
\begin{equation}
h^{0}(K_{M}+A)\leq 2. \label{m2.3.8}
\end{equation}
We note that
\begin{equation}
1\leq (K_{M}+A)^{2}A. \label{m2.3.13}
\end{equation}
Since $g_{2}(M,A)\geq h^{1}(\mathcal{O}_{M})$ by Theorem \ref{Theorem1.4} and $g_{1}(M,K_{M}+A,A)=1+(K_{M}+A)^{2}A$, we have 
\begin{eqnarray}
&&h^{0}(2K_{M}+2A)-h^{0}(2K_{M}+A) \label{m2.3.14}\\
&&\geq g_{1}(M,K_{M}+A,A) \nonumber \\
&&=1+(K_{M}+A)^{2}A \nonumber
\end{eqnarray}
and 
\begin{equation}
(K_{M}+A)^{2}A\leq 3 \label{m2.3.15}
\end{equation}
by Proposition \ref{32P-T1} (\ref{m2.3.1}) since $h^{0}(2K_{M}+2A)=4$.
\\
\par
Here we divide the argument into three cases.
\\
(i) Assume that $(K_{M}+A)^{2}A=1$.
Then $(K_{M}+A)A^{2}=1$ and $A^{3}=1$ by Proposition \ref{GHIT}.
So we get $g(M,A)=2$ and this is the type (i) in Claim \ref{T-CL1}.
\\
(ii) Assume that $(K_{M}+A)^{2}A=2$.
Then $g_{1}(M,K_{M}+A,A)=3$.
By Proposition \ref{GHIT}, we have $1\leq (K_{M}+A)^{3}\leq 4$.
Hence by Proposition \ref{32P-T1} (\ref{m2.3.1}) and Theorem \ref{Theorem1.4} we have 
\begin{equation}
d_{1}=0, 1. \label{m2.3.16}
\end{equation}
By (\ref{m2.3.8}), Proposition \ref{32P-T1} (3) and the assumption $h^{0}(2K_{M}+2A)=4$ we have
\begin{eqnarray*}
2&\leq& h^{0}(2(K_{M}+A))-h^{0}(K_{M}+A) \\
&=&d_{2}+g_{1}(M,K_{M}+A,A)\\
&=&d_{2}+3.
\end{eqnarray*}
Namely we have 
\begin{equation}
-1\leq d_{2}. \label{m2.3.17}
\end{equation}
By Proposition \ref{32P-T1} (\ref{m2.3.2}) and the assumption $h^{0}(2K_{M}+2A)=4$ we have
\begin{eqnarray*}
4&\geq& h^{0}(2(K_{M}+A))-h^{0}(K_{M}+A) \\
&=&d_{2}+3.
\end{eqnarray*}
Namely we have 
\begin{equation}
1\geq d_{2}. \label{m2.3.18}
\end{equation}
So we get the following table by the same argument as in the proof of Theorem \ref{EC1}.

\begin{center}
\begin{tabular}{ccccccc}\hline
 & $d_{1}$ & $d_{2}$ & $d_{2}-d_{1}$ & $(K_{M}+A)^{2}K_{M}$ & $(K_{M}+A)^{3}$ & $\chi(\mathcal{O}_{M})$\\
\hline
(2.1) & $0$ & $-1$ & $-1$ & $4\chi(\mathcal{O}_{M})-2$ & $4\chi(\mathcal{O}_{M})$ & $1$ \\
(2.2) & $0$ & $0$ & $0$ & $4\chi(\mathcal{O}_{M})$ & $4\chi(\mathcal{O}_{M})+2$ & $0$ \\
(2.3) & $0$ & $1$ & $1$ & $4\chi(\mathcal{O}_{M})+2$ & $4\chi(\mathcal{O}_{M})+4$ & $0$ \\
(2.4) & $1$ & $-1$ & $-2$ & $4\chi(\mathcal{O}_{M})-4$ & $4\chi(\mathcal{O}_{M})-2$ & $1$ \\
(2.5) & $1$ & $0$ & $-1$ & $4\chi(\mathcal{O}_{M})-2$ & $4\chi(\mathcal{O}_{M})$ & $1$ \\
(2.6) & $1$ & $1$ & $0$ & $4\chi(\mathcal{O}_{M})$ & $4\chi(\mathcal{O}_{M})+2$ & $0$ \\
\end{tabular}
\end{center}
\noindent
(ii.1) First we consider the case (2.4).
Then $(K_{M}+A)^{3}=2$.
By Proposition \ref{GHIT} we have 
\begin{eqnarray*}
4&=&((K_{M}+A)^{2}A)^{2} \\
&\geq&((K_{M}+A)^{3})((K_{M}+A)A^{2}) \\
&=&2(K_{M}+A)A^{2}.
\end{eqnarray*}
Hence
$(K_{M}+A)A^{2}\leq 2$.
\\
(ii.1.1) If $(K_{M}+A)A^{2}=2$, then we also see that 
\begin{eqnarray*}
4&\geq&((K_{M}+A)A^{2})^{2} \\
&\geq&(A^{3})((K_{M}+A)^{2}A) \\
&=&2A^{3}.
\end{eqnarray*}
Therefore
$A^{3}\leq 2$.
But since $(K_{M}+2A)A^{2}$ is even and $A^{3}>0$, we have $A^{3}=2$.
Hence $(A^{3})((K_{M}+A)^{2}A)=((K_{M}+A)A^{2})^{2}$.
By \cite[Corollary 2.5.4]{BeSo-Book} we have $K_{M}+A\equiv A$, that is, $K_{M}\equiv 0$.
In particular, $g_{2}(M,A)=g_{2}(M,K_{M}+A)$.
But since $d_{1}\neq d_{2}$ in the case (2.4), this is impossible.
\\
\\
(ii.1.2) If $(K_{M}+A)A^{2}=1$, then $A^{3}=1$ by Proposition \ref{GHIT}.
Hence we see that $g(M,A)=2$ and this is the type (i) in Claim \ref{T-CL1}.
\\
\\
(ii.2) Next we consider the cases (2.1), (2.3) and (2.5).
Then $(K_{M}+A)^{3}=4$.
Since $(K_{M}+A)^{2}A=2$, by Proposition \ref{GHIT}, we have $(K_{M}+A)A^{2}=1$ and by Proposition \ref{GHIT} we have $1=((K_{M}+A)A^{2})^{2}\geq ((K_{M}+A)^{2}A)(A^{3})\geq 2A^{3}$.
Since $A^{3}>0$, this is impossible.
\\
\\
(ii.3) Next we consider the cases (2.2) and (2.6).
Then $(K_{M}+A)^{3}=2$.
By Proposition \ref{GHIT}, we have $(K_{M}+A)A^{2}\leq 2$ since $(K_{M}+A)^{2}A=2$.
\\
\\
(ii.3.1) If $(K_{M}+A)A^{2}=2$, then by the same argument as (ii.1.1) above, we have $K_{M}\equiv 0$.
In this case 
$$h^{0}(2K_{M}+A)=\chi(2K_{M}+A)=\chi(A)=\frac{1}{6}A^{3}+\frac{1}{12}c_{2}(M)A$$and
$$h^{0}(2K_{M}+2A)=\chi(2K_{M}+2A)=\chi(2A)=\frac{4}{3}A^{3}+\frac{1}{6}c_{2}(M)A.$$
Since $(K_{M}+A)^{3}=2$ and $K_{M}\equiv 0$, we have $A^{3}=2$
and $g_{1}(M,K_{M}+A,A)=g(M,A)=3$.
By Proposition \ref{32P-T1} (\ref{m2.3.1}) we have 
\begin{eqnarray*}
h^{0}(2K_{M}+A)
&=&h^{0}(2K_{M}+2A)-d_{1}-g_{1}(M,K_{M}+A,A)\\
&=&1-d_{1}.
\end{eqnarray*}
Hence we get $d_{1}=0, 1$ because $d_{1}\geq 0$ 
by Theorem \ref{Theorem1.4}. 
\\
\\
(ii.3.1.1) If $d_{1}=1$, then $h^{0}(2K_{M}+A)=0$ and $(1/6)A^{3}+(1/12)c_{2}(M)A=0$.
Therefore $h^{0}(2K_{M}+2A)=A^{3}=2$ and this is impossible.
\\
\\
(ii.3.1.2) If $d_{1}=0$, then $h^{0}(2K_{M}+A)=1$ and $(1/6)A^{3}+(1/12)c_{2}(M)A=1$.
Hence $h^{0}(2K_{M}+2A)=A^{3}+2=4$.
We note that $h^{i}(A)=0$ for every positive integer $i$ 
because $K_{M}\equiv 0$.
Hence $1=h^{0}(2K_{M}+A)=\chi(2K_{M}+A)=\chi(A)=h^{0}(A)$.
So this is the type (ii) in Claim \ref{T-CL1}.
\\
\\
(ii.3.2) If $(K_{M}+A)A^{2}=1$, then 
\begin{eqnarray*}
1&=&((K_{M}+A)A^{2})^{2} \\
&\geq&((K_{M}+A)^{2}A)(A^{3}) \\
&=&2A^{3}
\end{eqnarray*}
and this is impossible.
\\
\\
(iii) Assume that $(K_{M}+A)^{2}A=3$.
Then $(K_{M}+A)^{3}\leq 9$ by Proposition \ref{GHIT} and $g_{1}(M,K_{M}+A,A)=1+(K_{M}+A)^{2}A=4$.
Since $h^{0}(2K_{M}+2A)=4$ in this case, we have $d_{1}=0$ by Proposition \ref{32P-T1} (\ref{m2.3.1}) and Theorem \ref{Theorem1.4}.
Moreover we see that $-2\leq d_{2}\leq 0$ by (\ref{m2.3.8}) and Proposition \ref{32P-T1} (\ref{m2.3.2}).
Since $d_{2}-d_{1}+2\chi(\mathcal{O}_{M})=(1/2)(K_{M}+A)^{2}K_{M}$ (see (\ref{m2.3.7})), we have

\begin{center}
\begin{tabular}{cccccc}\hline
 & $d_{1}$ & $d_{2}$ & $d_{2}-d_{1}$ & $(K_{M}+A)^{2}K_{M}$ & $(K_{M}+A)^{3}$ \\
\hline
(3.1) & $0$ & $-2$ & $-2$ & $4\chi(\mathcal{O}_{M})-4$ & $4\chi(\mathcal{O}_{M})-1$ \\
(3.2) & $0$ & $-1$ & $-1$ & $4\chi(\mathcal{O}_{M})-2$ & $4\chi(\mathcal{O}_{M})+1$ \\
(3.3) & $0$ & $0$ & $0$ & $4\chi(\mathcal{O}_{M})$ & $4\chi(\mathcal{O}_{M})+3$ \\
\end{tabular}
\end{center}

First we consider the case (3.1).
Since $1\leq (K_{M}+A)^{3}\leq 9$, we have
$(\chi(\mathcal{O}_{M}), (K_{M}+A)^{3})=(1,3)$ or $(2,7)$.
\par
Next we consider the case (3.2).
Then we have $(\chi(\mathcal{O}_{M}),(K_{M}+A)^{3})=(0,1)$, $(1,5)$ or $(2,9)$.
\\
Finally we consider the case (3.3).
In this case, we get $(\chi(\mathcal{O}_{M}),(K_{M}+A)^{3})=(0,3)$ or $(1,7)$.
\\
(iii.1) Here we note that
if $(K_{M}+A)^{3}\geq 5$, 
then by Proposition \ref{GHIT}
\begin{eqnarray*}
9&=&((K_{M}+A)^{2}A)^{2} \\
&\geq&((K_{M}+A)^{3})((K_{M}+A)A^{2}) \\
&\geq&5(K_{M}+A)A^{2}.
\end{eqnarray*}
and we have 
$(K_{M}+A)A^{2}=1$ and $A^{3}=1$ by Proposition \ref{GHIT}.
Hence $g(M,A)=2$ and this is the type (i) in Claim \ref{T-CL1}.
\\
\\
(iii.2) Next we consider the case where $(K_{M}+A)^{3}=3$.
By Proposition \ref{GHIT}, we see that $(K_{M}+A)A^{2}\leq 3$.
\par
If $(K_{M}+A)A^{2}\leq 2$, then $A^{3}=1$ because $(K_{M}+A)^{2}A=3$.
But since $(K_{M}+2A)A^{2}$ is even, we see that $(K_{M}+A)A^{2}=1$ and $A^{3}=1$. Namely we have $g(M,A)=2$ and this is the type (i) in Claim \ref{T-CL1}.
\par
So we may assume that $(K_{M}+A)A^{2}=3$.
Then $((K_{M}+A)A^{2})((K_{M}+A)^{3})=((K_{M}+A)^{2}A)^{2}=9$.
Here we will prove the following lemma.
\begin{Lemma}\label{32L-T1}
Let $X$ be a smooth projective variety of dimension $3$.
Let $D_{1}$, $D_{2}$ and $D_{3}$ be divisors on $X$.
Assume the following:
\\
{\rm (1)} $D_{1}^{2}D_{3}>0$.
\\
{\rm (2)} $D_{3}$ is semiample and big.
\\
{\rm (3)} $(D_{1}^{2}D_{3})(D_{2}^{2}D_{3})=(D_{1}D_{2}D_{3})^{2}$.
\\
{\rm (4)} $D_{1}^{2}D_{3}=D_{2}^{2}D_{3}$.
\\
Then $(D_{1}-D_{2})D_{3}D=0$ holds for any divisor $D$ on $X$.
\end{Lemma}
\noindent
{\em Proof.}
By the assumption (2), there exists a smooth surface $S\in |mD_{3}|$
for some $m>0$.
Then by the assumption (3) we have $(D_{1}|_{S})^{2}(D_{2}|_{S})^{2}=((D_{1}|_{S})(D_{2}|_{S}))^{2}$.
So by the assumptions (1) and (4) we have $D_{1}|_{S}\equiv D_{2}|_{S}$.
In particular $(D_{1}|_{S})(D|_{S})=(D_{2}|_{S})(D|_{S})$ for any divisor $D$ on $X$.
Therefore $D_{1}D(mD_{3})=D_{2}D(mD_{3})$.
Hence we get the assertion. $\Box$
\\
\par
Since $K_{M}+A$ is semiample and big, we see that
$(K_{M}+A)^{2}D=A(K_{M}+A)D$ for any divisor $D$ on $M$ by Lemma \ref{32L-T1}.
Therefore $K_{M}D(K_{M}+A)=0$ for any divisor $D$ on $X$.
\par
Next we calculate $h^{0}(2K_{M}+2A)$ and $h^{0}(K_{M}+A)$.
Then by the Hirzebruch-Riemann-Roch theorem and the Kodaira vanishing theorem we have
$$h^{0}(2K_{M}+2A)=4+(1/6)c_{2}(M)A-3\chi(\mathcal{O}_{M}),$$
and
$$h^{0}(K_{M}+A)=(1/2)+(1/12)c_{2}(M)A-\chi(\mathcal{O}_{M}).$$

Since we are considering the case where $(K_{M}+A)^{3}=3$, we have $\chi(\mathcal{O}_{M})=0$ or $1$.
\par
If $\chi(\mathcal{O}_{M})=0$, then $4=h^{0}(2K_{M}+2A)=4+(1/6)c_{2}(M)A$.
Hence $c_{2}(M)A=0$.
But then $h^{0}(K_{M}+A)=1/2$ and this is impossible.
\par
If $\chi(\mathcal{O}_{M})=1$, then $(M,A)$ satisfies the case (3.1) and 
$4=h^{0}(2K_{M}+2A)=4+(1/6)c_{2}(M)A-3\chi(\mathcal{O}_{M})=1+(1/6)c_{2}(M)A$.
Hence $c_{2}(M)A=18$ and $h^{0}(K_{M}+A)=1$.
On the other hand, by Theorem \ref{Theorem1.3} we have $1=h^{0}(K_{M}+A)=g_{2}(M,A)-h^{2}(\mathcal{O}_{M})+h^{3}(\mathcal{O}_{M})$.
Hence $g_{2}(M,A)=1+h^{2}(\mathcal{O}_{M})-h^{3}(\mathcal{O}_{M})$ and $d_{1}=\chi(\mathcal{O}_{M})=1$.
But $d_{1}=0$ in this case (3.1).
Hence this is also impossible.
\\
\\
(iii.3) Next we consider the case where $(K_{M}+A)^{3}=1$.
Then $(M,A)$ satisfies the case (3.2).
In particular $g_{2}(M,A)=h^{1}(\mathcal{O}_{M})$.
We also get $(K_{M}+A)^{2}K_{M}=-2$ from the assumption that $(K_{M}+A)^{2}A=3$ or $\chi(\mathcal{O}_{M})=0$.
In particular $\kappa(M)=-\infty$ and $h^{3}(\mathcal{O}_{M})=0$.
Here we note $g_{1}(M,K_{M}+A)=1+(1/2)(3K_{M}+2A)(K_{M}+A)^{2}=1$.
We also note that $h^{1}(\mathcal{O}_{M})>0$ because $\kappa(M)=-\infty$ and $\chi(\mathcal{O}_{M})=0$.
Hence by \cite[(4.9) Corollary]{Fujita89} we have $h^{1}(\mathcal{O}_{M})=1$ and $(M,K_{M}+A)$ is birationally equivalent to $(V,H)$ which is a scroll over an elliptic curve because $K_{M}+A$ is nef and big.
This is the type (iii) in Claim \ref{T-CL1}.
\par
These complete the proof of Claim \ref{T-CL1}. $\Box$
\\
\par
Here we consider the case where $g(M,A)=2$.
In this case by the classification of $(M,A)$ with $g(M,A)=2$ (\cite[(1.10) Theorem and Section 2]{Fujita87-2}) we see that $(M,A)$ is one of the following type: $\mathcal{O}(K_{M})=\mathcal{O}_{M}$, $h^{1}(\mathcal{O}_{M})=0$, $h^{0}(A)>0$ and $A^{3}=1$.
\par
Then $h^{0}(A)=(1/6)A^{3}+(1/12)c_{2}(M)A$ and $h^{0}(2A)=(4/3)A^{3}+(1/6)c_{2}(M)A$.
Since $4=h^{0}(2K_{M}+2A)=h^{0}(2A)$, we have
$4=(4/3)A^{3}+(1/6)c_{2}(M)A=(4/3)+(1/6)c_{2}(M)A$.
Hence $c_{2}(M)A=16$.
But then $h^{0}(A)=3/2$ and this is impossible.
\par
Therefore $(M,A)$ is one of the types (1) and (2) in Theorem \ref{EC2}.
\\
\\
($\beta$) Assume that $(M,A)$ satisfies one of the types (1) and (2) in Theorem \ref{EC2}.
\\
($\beta$.1) Assume that $(M,A)$ satisfies the type (1) in Theorem \ref{EC2}.
Here we note that $h^{i}(A)=0$ for every positive integer $i$.
Then
\begin{eqnarray*}
h^{0}(A)
&=&\chi(A)\\
&=&\frac{1}{6}A^{3}+\frac{1}{12}c_{2}(M)A.
\end{eqnarray*}
Hence we have $c_{2}(M)A=8$.
Therefore 
\begin{eqnarray*}
h^{0}(2K_{M}+2A)
&=&\chi(2K_{M}+2A)\\
&=&\chi(2A)\\
&=&\frac{4}{3}A^{3}+\frac{1}{6}c_{2}(M)A\\
&=&4.
\end{eqnarray*}
\noindent
\\
($\beta$.2) Assume that $(M,A)$ satisfies the type (2) in Theorem \ref{EC2}.
\par
First we note the following.
\begin{Claim}
$h^{0}(2K_{M}+A)=0$.
\end{Claim}
\noindent
{\em Proof.}
Since $(M,K_{M}+A)$ is birationally equivalent to a scroll $(V,H)$ over a smooth ellitpic curve $B$,
there exist a smooth projective $3$-fold $T$ and birational morphisms $\mu: T\to M$ and $\nu: T\to V$ such that $\mu^{*}(K_{M}+A)=\nu^{*}(H)$.
Here we note that $V$ is smooth.
Then $h^{0}(2K_{M}+A)=h^{0}(\mu^{*}(2K_{M}+A))=h^{0}(K_{T}+\mu^{*}(K_{M}+A))=h^{0}(K_{T}+\nu^{*}(H))=h^{0}(\nu^{*}(K_{V}+H))=h^{0}(K_{V}+H)=0$.
This completes the proof. $\Box$
\\
\par
We also see that
$g_{1}(M,K_{M}+A,A)=1+(K_{M}+A)^{2}A=4$.
Hence from Proposition \ref{32P-T1} (\ref{m2.3.1}) we get
\begin{eqnarray*}
h^{0}(2(K_{M}+A))
&=&h^{0}(2K_{M}+A)+g_{2}(M,A)-h^{1}(\mathcal{O}_{M})+g_{1}(M,K_{M}+A,A)\\
&=&4.
\end{eqnarray*}
Therefore we get the assertion of Theorem \ref{EC2}. $\Box$

\begin{Remark}
By Theorem \ref{EC2}, we see that if $\kappa(K_{X}+L)=3$ and $h^{0}(2(K_{X}+L))=4$, then $h^{0}(K_{X}+L)=1$.
\end{Remark}

\begin{Example}
Here we give an example of this case.
\begin{itemize}
\item [\rm (1)] An example of the type (1) in Theorem \ref{EC2}.
In \cite[Theorem 1.1]{Beauville}, Beauville gave an example of
a polarized Calabi-Yau threefold $(X,L)$ such that 
$h^{0}(L)=1$ and $L^{3}=2$. This is an example.
For details, see \cite[Theorem 1.1]{Beauville}.
\item [\rm (2)] An example of the type (2) in Theorem \ref{EC2}.
Let $C$ be an elliptic curve and let $\mathcal{E}$ be an ample vector bundle of rank $3$ on $C$ with $c_{1}(\mathcal{E})=1$.
Then $\mathcal{E}$ is indecomposable.
We note that such a vector bundle exists.
Let $M=\mathbb{P}_{C}(\mathcal{E})$ and $A=4H(\mathcal{E})-f^{*}(c_{1}(\mathcal{E}))$, where $f:M\to C$ is the natural map.
Then by \cite[Theorem 3.1]{Miyaoka87} we see that $A$ is ample,
and we also see that $(M,K_{M}+A)$ is a scroll over a smooth elliptic curve.
We can also check that $h^{0}(K_{M}+A)=h^{0}(H(\mathcal{E}))=1$, $h^{2}(\mathcal{O}_{M})=0$, $h^{1}(\mathcal{O}_{M})=1$, $g_{2}(M,A)=1$, $g_{1}(M,K_{M}+A,A)=4$ and $h^{0}(2K_{M}+A)=0$.
Therefore by Proposition \ref{32P-T1} (2) we have
$h^{0}(2K_{M}+2A)=h^{0}(2K_{M}+A)+g_{2}(M,A)-h^{1}(\mathcal{O}_{M})+g_{1}(M,K_{M}+A,A)=4$.
\end{itemize}
\end{Example}

\begin{flushright}
Yoshiaki Fukuma \\
Department of Mathematics \\
Faculty of Science \\
Kochi University \\
Akebono-cho, Kochi 780-8520 \\
Japan \\
E-mail: fukuma@kochi-u.ac.jp
\end{flushright}

\end{document}